\RequirePackage{etoolbox}
\csdef{input@path}{{style/}{graphics/}}
\documentclass[MSNbibl,number,citesort,dvips]{arxbj}
\usepackage{mathbh}
\usepackage{upgreek}
\usepackage{graphicx}

\aid{0}
\volume{21}
\issue{1}
\pubyear{2015}
\firstpage{276}
\lastpage{302}
\doi{10.3150/13-BEJ567} % kopijuoti is 'New paper accepted'

\makeatletter

\newcommand{\rrvert}{\vert}
\newcommand{\llvert}{\vert}
\newcommand{\cf}{cf.\ }

\newcommand{\resp}{resp.\ }

\newcommand{\iid}{i.i.d.\ }
\newcommand{\fdd}{f.d.d.\ }
\newcommand{\wrt}{w.r.t.\ }
\newcommand{\rhs}{r.h.s.\ }
\newcommand{\RR}{\mathbb{R}}
\newcommand{\NN}{\mathbb{N}}
\newcommand{\PP}{\mathbb{P}}
\newcommand{\EE}{\mathbb{E}}
\newcommand{\FF}{\mathbb{F}}
\newcommand{\D}{\mathrm{d}}
\newcommand{\erf}{\operatorname{erf}}
\newcommand{\pr}{\operatorname{pr}}
\newcommand{\Borel}{{\mathcal B}}
\newcommand{\sigmaA}{{\mathcal A}}
\newcommand{\Half}{{\mathcal H}}
\newcommand{\Plane}{{\mathcal E}}
\newcommand{\DepSet}{{\mathcal K}}
\newcommand{\PotentialDepSet}{{\mathcal L}}
\newcommand{\finite}{\mathcal{F}}
\newcommand{\powerset}{\mathcal{P}}
\newcommand{\eins}{\mathbf{1}}

\newcommand{\pmid}{\dvtx }
\newtheorem{theorem}{Theorem}

\newtheorem{corollary}[theorem]{Corollary}
\newtheorem{lemma}[theorem]{Lemma}

\newproclaim{definition}[theorem]{Definition}
\newremark{example}[theorem]{Example}

\newremark{remark}[theorem]{Remark}
\makeatother

\begin{document}
\begin{frontmatter}

\title{An exceptional max-stable process fully parameterized by its
extremal coefficients}
\runtitle{A max-stable process parameterized by its ECF}

\begin{aug}
%%%% inicialai - be tarpu
\author[1]{\inits{K.}\fnms{Kirstin} \snm{Strokorb}\corref{}\thanksref{e1}\ead[label=e1,mark]{strokorb@math.uni-mannheim.de}} \and
\author[1]{\inits{M.}\fnms{Martin} \snm{Schlather}\thanksref{e2}\ead[label=e2,mark]{schlather@math.uni-mannheim.de}}
%%\runauthor{} %% auto
\address[1]{Institute of Mathematics, University of Mannheim, 68131
Mannheim, Germany.\\ \printead{e1,e2}}
\end{aug}

% HISTORY:
\received{\smonth{3} \syear{2013}}
\revised{\smonth{9} \syear{2013}}

% ABSTRACT
%
\begin{abstract}
The extremal coefficient function (ECF) of
a max-stable process $X$ on some index set $T$
assigns to each finite subset $A \subset T$ the effective
number of independent random variables among the
collection $\{X_t\}_{t \in A}$.
We introduce the class of Tawn--Molchanov processes
that is in a 1:1 correspondence with the class of ECFs,
thus also proving a complete characterization of the ECF
in terms of negative definiteness. The corresponding
Tawn--Molchanov process turns out to be exceptional among
all max-stable processes sharing the same ECF in that its
dependency set is maximal \wrt inclusion. This entails sharp
lower bounds for the finite dimensional distributions of
arbitrary max-stable processes in terms of its ECF.
A spectral representation of the Tawn--Molchanov process
and stochastic continuity are discussed. We also show
how to build new valid ECFs from given ECFs by means of
Bernstein functions.
\end{abstract}

% KEYWORDS
% visi is mazosios raides ir pagal abecele
%

\end{frontmatter}

%s1 #&#
\section{Introduction}\label{sect:intro}

Besides the class of {square integrable processes}, the class of
temporal or spatial \emph{max-stable processes} is of common interest
in stochastics and statistics, \cf\cite
{dehaan_84,ginehahnvatan_90,wangstoev_10,blanchetdavison_11,buishandetalii_08,naveauetalii_09},
for example. In spite of considerable differences between these two
classes, for example, the non-existence of the first moments in case of
max-stable processes with \emph{unit Fr\'echet marginals}, connections
between the two classes have been made for instance, by the \emph
{extremal Gaussian process} \cite{schlather_02} and the \emph
{Brown--Resnick process} \cite{kabluchkoetalii_09} that are
parameterized by a {correlation function} and a {variogram}, respectively.

Naturally, extremal dependence measures such as the \emph{extremal
coefficients} \cite{smith_90,schlathertawn_02}, the \emph{(upper) tail
dependence coefficients} \cite
{beirlantetalii_03,davismikosch_09,falk_05,colesheffernantawn_99} or
other special cases of the \emph{extremogram} \cite{davismikosch_09}
are appropriate summary statistics for max-stable processes.
In this article, we capture the full set of extremal coefficients of a
max-stable process $X=\{X_t\}_{t \in T}$ on some space $T$ in the
so-called \emph{extremal coefficient function (ECF)} $\theta$, which
assigns to each finite subset $A$ of $T$ the effective number of
independent variables among the collection $\{X_t\}_{t \in A}$. We
introduce a subclass of max-stable processes that is parameterized by
the ECF, and thus reveal some analogies to \emph{Gaussian processes}
and \emph{positive definite functions} as follows:

Among (zero mean) square integrable processes, the subclass of \emph
{Gaussian processes} takes a unique role, since it is in a 1--1
correspondence with the set of \emph{covariance functions}, which are
precisely the \emph{positive definite functions}. This fact can be
proven by means of Kolmogorov's extension theorem and is illustrated in
the following graph:%\vspace*{6pt}

%
%f1 #&#
\begin{center}
\begin{tabular}{c}

\includegraphics{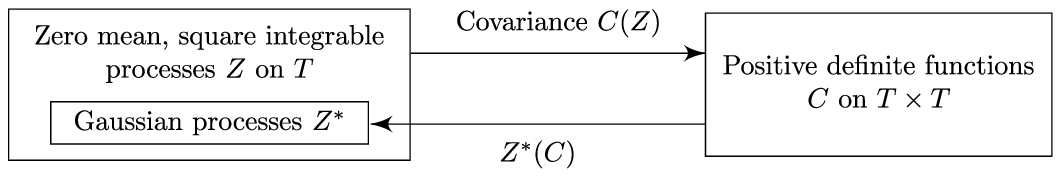}

\end{tabular}\vspace*{6pt}
\end{center}
%{\begin{minipage}{4.7cm}
%Zero mean, square integrable \\
%processes $Z$ on $T$ \\[2mm]
%{\begin{minipage}{4cm}
%$\phantom{\tiny{a}}$\\
%Positive definite functions\\
%$C$ on $T \times T$\\
%$\phantom{\tiny{a}}$
%$C(Z)$} +(3.75,0.4);
%+(-4.25,-0.5);

\noindent In case $T$ is a metric space, the Gaussian process $Z^*(C)$ is
continuous in the mean square sense (and then also stochastically
continuous) if and only if the covariance function $C$ is continuous if
and only if $C$ is continuous on the diagonal (\cf\cite{scheuerer_10}, Theorem
5.3.3). Well-known operations on the set of positive
definite functions $C$, and hence on the corresponding Gaussian
processes $Z^*(C)$, include convex combinations and pointwise limits.
Moreover, \emph{Bernstein functions} play an important role for the
construction of positive definite functions.

In our case, the crucial role of zero mean Gaussian processes is taken
by the class of \emph{Tawn--Molchanov processes (TM processes)}, which
are in fact the spatial generalization of the multivariate \emph
{max-linear} model of \cite{schlathertawn_02}. Using Kolmogorov's
extension theorem, we shall see that each ECF $\theta$ (of some
max-stable process) uniquely determines a TM process $X^*(\theta)$
having the same ECF (Theorem \ref{thm:ECF_ND}). Alongside, we
generalize a multivariate result
\cite{molchanov_08}, Corollary 1, to the spatial setting, proving that
the ECFs coincide with the functions $\theta$ on $\finite(T)$ (the
\emph
{set of finite subsets} of $T$) that are normalized to $\theta
(\varnothing
)=0$ and $\theta(\{t\})=1$ for $t \in T$ and that are \emph{negative
definite} (or equivalently \emph{completely alternating}) in a sense to
be explained below (\cf Definition \ref{def:ND_CA}). This can be
illustrated in analogy to the above sketch:\vspace*{6pt}

%f2 #&#
\begin{center}
\begin{tabular}{c}

\includegraphics{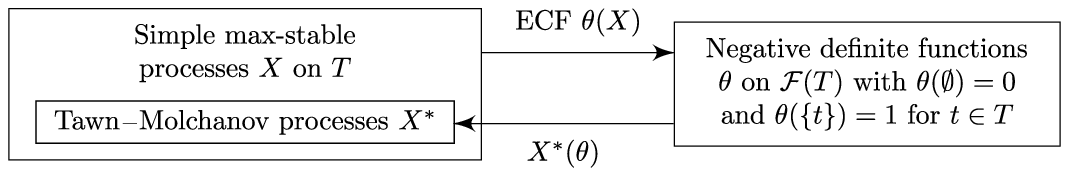}

\end{tabular}\vspace*{6pt}
\end{center}
%
%{\begin{minipage}{5.6cm}
%Simple max-stable\\
%processes $X$ on $T$ \\[2mm]
%{\begin{minipage}{4.5cm}
%Negative definite functions\\
%$\theta$ on $\finite(T)$ with $\theta(\varnothing)=0$\\
%and $\theta(\{t\})=1$ for $t \in T$
%+(2.45,0.4);
%+(-2.8,-0.5);

\noindent Having identified the ECF $\theta$ as a negative definite quantity
allows for several immediate consequences:
First, we obtain an \emph{integral representation} of $\theta$ as a
mixture of maps $A \mapsto\mathbh{1}_{A \cap Q \neq\varnothing}$
(Corollary \ref{cor:intrep}) and derive a \emph{spectral
representation} for the corresponding TM process $X^*(\theta)$ (Theorem
\ref{thm:spectralrep}).
Second, we consider operations on ECFs that allow to build new ECFs
from given ones. We find that ECFs allow for convex combinations and
pointwise limits (Corollaries \ref{cor:Theta_convex} and \ref
{cor:Theta_compact}) and that the class of \emph{Bernstein functions}
operates on ECFs (Corollary \ref{cor:opBernstein}). We also recover the
``triangle inequalities'' for $\theta$ from \cite{cooleyetalii_06}, Proposition
4,  and see that the inequalities therein correspond to
three specific choices of a Bernstein function, whereas we may plug in
arbitrary Bernstein functions to obtain further ``triangle
inequalities'' (Corollary \ref{cor:triangleineq}).

For $T$ being a metric space, we discuss \emph{stochastic continuity}:
The TM process $X^*(\theta)$ is stochastically continuous if and only
if $\theta$ is continuous (\cf Definition \ref{def:ECF_cts}) if and
only if the bivariate map $(s,t) \mapsto\theta(\{s,t\})$ is continuous
if and only if the bivariate map $(s,t) \mapsto\theta(\{s,t\})$ is
continuous on the diagonal (Theorem \ref{thm:ECprocess_cty}).

\begin{figure}

\includegraphics{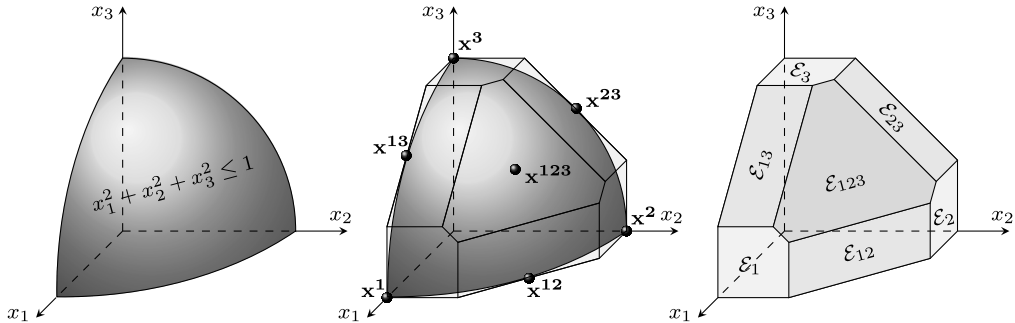}

\caption{Examples of dependency sets in a trivariate setting:
a ``typical'' dependency set $\DepSet$ (left) and a dependency set
$\DepSet^*$ stemming from a TM process (right). It is shown that
$\DepSet\subset\DepSet^*$ (middle). For further details, see
the introduction, Example \protect\ref{example:ball}, Lemma
\protect\ref{lemma:DepSet} and Theorem \protect\ref{thm:starDepSet}.}
\label{fig:ballDepSet}
\end{figure}

Finally, we address the exceptional role of the TM processes among
simple max-stable processes. To this end, Molchanov's \emph{dependency
set} $\DepSet$ \cite{molchanov_08} is transferred to max-stable
processes $X$. It comprises the finite dimensional distributions
(f.d.d.) of $X$ (Lemma \ref{lemma:DepSetprops}). Now, let $\DepSet
^*(\theta)$ denote the dependency set of the process $X^*(\theta)$.
Then we identify $\DepSet^*(\theta)$ as intersection of halfspaces that
are directly given by the ECF $\theta$ (Theorem \ref{thm:starDepSet}).
It turns out that $\DepSet^*(\theta)$ is exceptional among the
dependency sets $\DepSet$ of all max-stable processes sharing the same
ECF $\theta$, since $\DepSet^*(\theta)$ is maximal \wrt inclusion as
illustrated in Figure~\ref{fig:ballDepSet}.
Since inclusion of dependency sets corresponds to stochastic ordering,
this observation leads to sharp inequalities for the \fdd of max-stable
processes in terms of its ECF $\theta$ (Corollary \ref{cor:fddinequalities}).

The text is structured as follows. After the introductory Section~\ref{sect:foundations}, the characterization of ECFs and the existence of
TM processes is established in Section~\ref{sect:TMprocess}. Section~\ref{sect:consequences} collects several immediate consequences and
related results, while Section~\ref{sect:depset} exhibits the
exceptional role of TM processes. Sections~\ref{sect:consequences} and
\ref{sect:depset} can be read independently.

%f3 #&#

%s2 #&#
\section{Foundations and definitions}\label{sect:foundations}

%s2.1 #&#
\subsection{Notation for max-stable processes and ECFs}

A stochastic process $X=\{X_t\}_{t \in T}$ on an arbitrary index set
$T$ is said to be \emph{max-stable} if for each $n \in\NN$ and
independent copies $X^{(1)},\ldots,X^{(n)}$ of $X$ the process of the
maxima $\{\bigvee_{i=1}^n X^{(i)}\}_{t \in T}$ has the same law as $\{
a_n(t) X_t + b_n(t)\}_{t \in T}$ for suitable norming functions
$a_n(t)>0$ and $b_n(t)$ on $T$. Without loss of generality, we shall
deal with max-stable processes that have \emph{standard Fr{\'e}chet
marginals}, that is, $\PP(X_t \leq x)=\mathrm{e}^{-1/x}$ for $t \in T$ and $x
\geq0$, and set $a_n(t)=n$ and $b_n(t)=0$. Such processes are called
\emph{simple} max-stable processes.

It has been shown (\cf\cite{kabluchko_09,dehaan_84,stoev_08})
that (simple) max-stable processes $X=\{X_t\}_{t \in T}$
allow for a \emph{spectral representation} $(\Omega,\sigmaA,\nu,V)$:
there exists a (sufficiently rich) measure space $(\Omega,\sigmaA,\nu
)$
and measurable functions $V_t:\Omega\rightarrow\RR_+$ (with
$\int_{\Omega} V_t(\omega) \nu(\D\omega) = 1$ for each $t \in T$),
such that the law of $X=\{X_t\}_{t \in T}$ equals the law of
%
%e1 #&#
\begin{equation}
\label{eqn:spectralrepresentation} \biggl\{ \bigvee_{(U,\omega)\in\Pi} U
V_t(\omega) \biggr\}_{t \in T}.
\end{equation}
Here $\Pi$ denotes a Poisson point process on $\RR_+ \times\Omega$
with intensity $u^{-2}\, \D u \times\nu(\D\omega)$. The functions $\{
V_t\}_{t \in T}$ are called \emph{spectral functions} and the measure
$\nu$ is called \emph{spectral measure}.

In order to describe the \emph{finite dimensional distributions}
(f.d.d.) of $X$, we shall fix some convenient notation first: Let $M
\subset T$ be some non-empty finite subset of $T$. By $\RR^M$ (\resp
$[0,\infty]^M$) we denote the space of real-valued (\resp$[0,\infty
]$-valued) functions on $M$. Elements of these spaces are denoted by
$x=(x_t)_{t \in M}$ where $x_t=x(t)$.
The standard scalar product is given through $\langle x,y \rangle=
\sum_{t \in M} x_t   y_t$. For a subset $L \subset M$, we write $\eins_L$
for the indicator function of $L$ in $\RR^M$ (\resp$[0,\infty]^M$),
such that $\{\eins_{\{t\}}\}_{t \in M}$ forms an orthonormal basis of
$\RR^M$. In this sense, the origin of $\RR^M$ equals $\eins
_{\varnothing
}$ being zero everywhere on $M$.
Using this notation, we emphasize the fact that a multivariate
distribution of a stochastic process is not any $|M|$-variate
distribution, but it is bound to certain points (forming the set $M$)
in the space $T$.
Finally, we consider some \emph{reference norm} $\| \cdot\|$ on $\RR
^M$ (not necessarily the one from the scalar product) and denote the
positive unit sphere $S_M:=\{ a \in[0,\infty)^M \pmid\| a \| = 1 \}$.

In terms of a spectral representation $(\Omega,\sigmaA,\nu,V)$, the
\fdd of $X$ are given through
%
%e2 #&#
\begin{equation}
\label{eqn:fddspectralrep} -\log\PP(X_t\leq x_t, t \in M) = \int
_{\Omega} \biggl(\bigvee_{t
\in M}
\frac{V_{t}(\omega)}{x_t} \biggr) \nu(\D\omega)
\end{equation}
for $x \in[0,\infty)^M\setminus\{\eins_\varnothing\}$.
Alternatively, the \fdd of $X$ for a finite subset $\varnothing\neq M
\subset T$ may be described by means of one of the following three
quantities that are all equivalent to the knowledge of the f.d.d.:
\begin{itemize}
\item the \emph{(finite dimensional) spectral measure} $H_M$ (\cf
\cite
{dehaanresnick_77,resnick_08}), that is, the Radon measure on $S_M$
such that
for $x \in[0,\infty)^M\setminus\{\eins_\varnothing\}$
%
%e3 #&#
\begin{equation}
\label{eqn:fddspectralmeasure} - \log\PP(X_{t} \leq x_t, t \in M) = \int
_{S_M} \biggl(\bigvee_{t
\in
M}
\frac{a_t}{x_t} \biggr) H_M (\D a),
\end{equation}
\item the \emph{stable tail dependence function} $\ell_M$ (\cf\cite
{beirlantetalii_03}), that is, the function on $[0,\infty)^M$ defined through
%
%e4 #&#
\begin{equation}
\label{eqn:fddstabletaildepfn} \ell_M(x) := - \log\PP(X_{t}
\leq1/x_t, t \in M) = \int_{S_M} \biggl( \bigvee
_{t \in M} a_t x_t \biggr)
H_M (\D a),
\end{equation}
\item the \emph{(finite dimensional) dependency set} $\DepSet_M$ (\cf
\cite{molchanov_08}), that is, the largest compact convex set
$\DepSet
_M \subset[0,\infty)^M$ satisfying
%
%e5 #&#
\begin{equation}
\label{eqn:DepSetdefn} \ell_M(x)=\sup \bigl\{ \langle x,y \rangle\pmid y \in
\DepSet_M \bigr\}\quad\quad \forall x \in[0,\infty)^M.
\end{equation}
\end{itemize}
In order to be a valid finite dimensional spectral measure of a simple
max-stable random vector $\{X_t\}_{t \in M}$, the only constraint that
a Radon measure $H_M$ on $S_M$ has to satisfy is that
\[
\int_{S_M} a_t H_M(\D a)=1
\]
for each $t \in M$. This ensures standard Fr{\'e}chet marginals.

Given a simple max-stable process $X$ on $T$, we may assign to a
non-empty finite subset $A \subset T$ the \emph{extremal coefficient}
$\theta(A)$ (\cf\cite{smith_90,schlathertawn_02}), that is
%
%e6 #&#
\begin{equation}
\label{eqn:ECFintro} \theta(A):= \lim_{x \to\infty}\frac{\log\PP (\bigvee_{t \in
A} X_t
\leq x )}{\log\PP(X_t \leq x)} = \int
_{S_M} \biggl(\bigvee_{t
\in A}
{a_t} \biggr) H_M (\D a) = \ell_M (
\eins_A ).
\end{equation}
Indeed, the expression $\log\PP(\bigvee_{t \in A} X_t \leq x)/\log\PP
(X_t \leq x)$ does not depend on $x$ and equals the right-hand side
(r.h.s.) for $A \subset M$. Observe that $\theta(A)$ takes values in
the interval $[1,|A|]$, where the value $1$ corresponds to full
dependence of the random variables $\{X_t\}_{t \in A}$ and the value
$|A|$ (number of elements in $A$) corresponds to independence.
Roughly speaking, the extremal coefficient $\theta(A)$ detects the
effective number of independent variables among the random variables $\{
X_t\}_{t \in A}$. It is coherent to set $\theta(\varnothing):=0$ to
obtain a function $\theta$ on $\finite(T)$, the \emph{set of finite
subsets of $T$}. We call the function
\[
\theta\dvtx \finite(T) \rightarrow[0,\infty)
\]
\emph{extremal coefficient function (ECF)} of the process $X$. The set
of all ECFs of simple max-stable processes on a set $T$ will be denoted by
%
%e7 #&#
\begin{equation}
\label{eqn:Theta_defn} \Theta(T)= \bigl\{\theta\dvtx \finite(T) \rightarrow[0, \infty)
\pmid \theta \mbox{ is an ECF of a simple max-stable process on } T  %
\bigr\}.
\end{equation}

%ex1 #&#
\begin{example}
The simplest ECFs are the functions $\theta(A)=|A|$ corresponding to a
process of independent random variables, and the indicator function
$\theta(A) = \mathbh{1}_{A \neq\varnothing}$ corresponding to a process
of identical random variables.
\end{example}

Rather sophisticated examples of ECFs can be obtained using spectral
representations $(\Omega,\sigmaA,\nu,V)$ of processes $X$. In these
terms the ECF $\theta$ of a process $X$ is given by
%
%e8 #&#
\begin{equation}
\label{eqn:ECFspectralrep} \theta(A)=\int_{\Omega} \biggl(\bigvee
_{t \in A} {V_t(\omega)} \biggr) \nu (\D\omega)
\end{equation}
for $A \in\finite(T) \setminus\{\varnothing\}$ and $\theta(\varnothing)=0$.

%ex2 #&#
\begin{example}[(Mixed Moving Maxima -- M3 process)]
Consider the simple max-stable stationary process $X$ on $\RR^d$ that
is given through the following spectral representation $(\Omega
,\sigmaA
,\nu,V)$:
\begin{itemize}
\item$(\Omega,\sigmaA,\nu)=(\FF\times\RR^d, {\mathcal F}\otimes
\Borel(\RR^d), \mu\otimes\D z)$, where $(\RR^d,\Borel(\RR^d),\D z)$
denotes the Lebesgue-measure on the Borel-$\sigma$-algebra of $\RR^d$
and where $(\FF,{\mathcal F},\mu)$ denotes a measure space of
non-negative measurable functions on $\RR^d$ with $\int_{\FF}
(\int_{\RR^d} f(z)\, \D z ) \mu(\D f) = 1$,
\item$V_t((f,z))=f(t-z)$ for $t \in\RR^d$,
\end{itemize}
then we call $X$ a \emph{Mixed Moving Maxima process (M3 process)}
(\cf
also \cite{kabluchkostoev_12,schlather_02,stoev_08,stoevtaqqu_06}).
Because of (\ref{eqn:ECFspectralrep}) the ECF $\theta$ of a Mixed
Moving Maxima process $X$ can be computed as
\[
\theta(A) =\int_{\FF} \int_{\RR^d} \biggl(
\bigvee_{t \in A} f(t-z) \biggr) \,\D z \mu (\D f)
\]
for $A \in\finite(\RR^d) \setminus\{\varnothing\}$ and $\theta
(\varnothing)=0$.
In case $\mu$ is a point mass at an indicator function~$f$ the
bivariate coefficient $\theta(\{s,t\})$ will be given by $\theta(\{
s,t\}
)=2-f*\check f(s-t)$, where $f*\check f$ means the convolution of $f$
with $\check f$ and $\check f(t)=f(-t)$.
\end{example}

%ex3 #&#
\begin{example}[(Brown--Resnick process)]\label{example:BR}
Consider the simple max-stable stationary process $X$ on $\RR^d$ that
is given through the following spectral representation $(\Omega
,\sigmaA
,\nu,V)$:
\begin{itemize}
\item$(\Omega,\sigmaA,\nu)$ denotes the probability space of a
Gaussian process $W$ on $\RR^d$ with stationary increments and
variogram $\gamma(z)=\EE(W_z-W_o)^2$ for $z \in\RR^d$.
\item$V_t(\omega)=\exp ({W_t(\omega)-\sigma^2(t)/2} )$
for $t
\in\RR^d$, where $\sigma^2(t)$ denotes the variance of $W_t$,
\end{itemize}
then we call $X$ a \emph{Brown--Resnick process} (\cf\cite
{kabluchkoetalii_09}). Because of (\ref{eqn:ECFspectralrep}) the ECF
$\theta$ of a Brown--Resnick process $X$ is
\[
\theta(A)= \EE\exp \biggl(\bigvee_{t \in A}
W_t- \sigma^2(t)/2 \biggr)
\]
for $A \in\finite(\RR^d) \setminus\{\varnothing\}$ and $\theta
(\varnothing)=0$.
Since the \fdd of $X$ only depend on the variogram $\gamma$, the
extremal coefficient $\theta(A)$ will also depend only on the values
$\{
\gamma(s-t)\}_{s,t \in A}$. In particular, we have $\theta(\{s,t\}
)=1+\erf(\sqrt{\gamma(s-t)/8})$ for the bivariate coefficient
$\theta(\{
s,t\})$, where
$\erf(x)=2/\sqrt{\uppi} \int_{0}^x \mathrm{e}^{-t^2} \,\mathrm{d}t$
denotes the error function (\cf\cite{kabluchkoetalii_09}). In case
the variogram equals $\gamma(z)=\lambda\| z \|^2$ for some $\lambda>
0$, explicit expressions for multivariate coefficients of higher orders
up to $d+1$ can be found in \cite{gentonetal_11}.
\end{example}

%s2.2 #&#
\subsection{A consistent max-linear model}

A multivariate simple max-stable distribution is called \emph
{max-linear} (or \emph{spectrally discrete}) if it arises as the
distribution of a random vector $X$ of the following form
\[
X_i = \bigvee_{j=1}^q
a_{ij} Z_j,\quad\quad i=1,\ldots,p,
\]
where $Z=\{Z_j\}^q_{j=1}$ is a vector of \iid unit Fr{\'e}chet random
variables and where $\{a_{ij}\}_{p \times q}$ is a matrix of
non-negative entries with $\sum_{j=1}^q a_{ij}=1$ for each row
$i=1,\ldots,p$. This is equivalent to requiring the spectral measure
$H_M$ from (\ref{eqn:fddspectralmeasure}) for $M=\{1,\ldots,\ldots,p\}$
to be the following \emph{discrete} measure on $S_M$
\[
H_M=\sum_{j=1}^q \|
a_j \| \delta_{a_j/\| a_j\|},
\]
where $a_j$ denote the column vectors of the matrix $\{a_{ij}\}_{p
\times q}$. Conversely, any discrete spectral measure of a simple
max-stable random vector gives rise to such a matrix. Surely, the ECF
of such a random vector $X=\{X_i\}_{i \in M}$ is
%
%e9 #&#
\begin{equation}
\label{eqn:ECFmaxlinear} \theta(A)=\sum_{j=1}^q
\bigvee_{i \in A}a_{ij}
\end{equation}
for $\varnothing\neq A \subset M$ and $\theta(\varnothing)=0$ (\cf(\ref
{eqn:ECFintro})).

In \cite{schlathertawn_02}, the authors introduce a max-linear model
for $X^*=\{X^*_i\}_{i \in M}$ where the column index $j$ ranges over
all non-empty subsets $L$ of $M$ and where non-negative coefficients
$\tau_L$ are given for each column $\varnothing\neq L \subset M$, more precisely
\[
X^*_i = \bigvee_{\varnothing\neq L \subset M} a_{i,L}
Z_L ,\quad\quad i \in M, \mbox{ with } a_{i,L} =
\tau_L \mathbh{1}_{i \in L},
\]
which is equivalent to
%
%e10 #&#
\begin{equation}
\label{eqn:taumodel} X^*_i = \bigvee_{i \in L \subset M}
\tau_{L} Z_L,\quad\quad i \in M.
\end{equation}
The\vspace*{2pt} model (\ref{eqn:taumodel}) is simple if and only if $\sum_{\varnothing\neq L \subset M} a_{iL} = \sum_{L \subset M \pmid i \in L}
\tau_L=1$ for each $i \in M$. It follows from (\ref{eqn:ECFmaxlinear})
that the ECF of model (\ref{eqn:taumodel}) is
\[
\theta(A)=\sum_{L \subset M \pmid A \cap L \neq\varnothing} \tau_L
\]
for $\varnothing\neq A \subset M$ and $\theta(\varnothing)=0$. Now, the
interesting aspect of this model (\ref{eqn:taumodel}) with given
coefficients $\tau_L$ is that such models are in 1--1 correspondence
with ECFs $\theta$ on the finite set $M$ (\cf\cite{schlathertawn_02}, Theorems 3 and
4). Alongside, this leads to a set of inequalities
which fully characterizes the set of ECFs $\Theta(M)$ for finite sets
$M$ (\cf\cite{schlathertawn_02}, Corollary 5). An alternative proof
for these inequalities is offered in \cite{molchanov_08}, Corollary 1,
and it is noticed therein that they are equivalent to a property called
complete alternation (see below).

\begin{figure}

\includegraphics{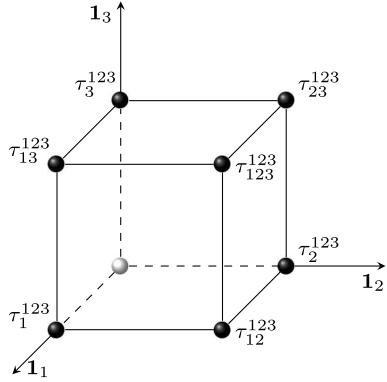}

\caption{Spectral measure representation of $\{X^*_{t}\}_{t \in M}$ for $M = \{
1,2,3\}$ if we choose the reference norm on $\RR^M$ to be the maximum
norm. In this case, the spectral measure simplifies to a sum of
weighted point masses on the vertices of a cube: $H^*_M = \sum_{\varnothing\neq L \subset M} \tau^M_L   \delta_{\eins_L}$.}

\label{fig:spectralmeasure}
\end{figure}

As we seek a spatial generalization of these results, let us consider a
max-stable processes $X^*=\{X^*_t\}_{t \in T}$ on an arbitrary index
set $T$, whose \fdd for a finite set $M$ are precisely of the above
form (\ref{eqn:taumodel}), where the coefficients $\tau_L$ now
additionally depend on $M$. That means we set the spectral measure
$H^*_M$ of the random vector $\{X^*_t\}_{t \in M}$
%
%e11 #&#
\begin{equation}
\label{eqn:starspectralmeasure} H^*_M := \sum_{\varnothing\neq L \subset M}
\tau^M_L \| \eins_L \| \delta_{\eins_L/\| \eins_L \|},
\end{equation}
such that the \fdd of the process $X^*$ are given by (\cf(\ref
{eqn:fddspectralmeasure}))
%
%e12 #&#
\begin{equation}
\label{eqn:starfdd} -\log\PP \bigl(X^*_t\leq x_t, t \in M
\bigr) =\sum_{\varnothing\neq L \subset M} \tau^M_L
\bigvee_{t \in L} \frac{1}{x_t}.
\end{equation}
Here $M$ ranges over all non-empty finite subsets of $T$, which we
express as $M \in\finite(T)\setminus\{\varnothing\}$.
Figure~\ref{fig:spectralmeasure} illustrates this spectral measure for
a trivariate distribution where $M=\{1,2,3\}$ in case the reference
norm is the maximum norm.

%f4 #&#

%le4 #&#
\begin{lemma} \label{lemma:tauprocess}
Let $T$ be an arbitrary set and let coefficients $\tau^M_L$ be given
for $M \in\finite(T)\setminus\{\varnothing\}$ and $L \in\finite
(M)\setminus\{\varnothing\}$, such that
\begin{enumerate}[(iii)]
\item[(i)]  $\tau^M_L \geq0$ for all $M \in\finite
(T)\setminus\{\varnothing\}$ and $L \in\finite(M)\setminus\{
\varnothing\}$,
\item[(ii)]  $\tau^M_L=\tau^{M \cup\{t\}}_L+\tau^{M
\cup\{
t\}}_{L \cup\{t\}}$ for all $M \in\finite(T)\setminus\{\varnothing\}$
and $L \in\finite(M)\setminus\{\varnothing\}$ and $t \in T\setminus M$,
\item[(iii)]  $\tau^{\{t\}}_{\{t\}}=1$ for all $t \in T$.
\end{enumerate}
Then the spectral measures $\{H^*_M\}_{M \in\finite(T) \setminus\{
\varnothing\}}$ from (\ref{eqn:starspectralmeasure}) define a simple
max-stable process $X^*=\{X^*_t\}_{t \in T}$ on $T$ with \fdd as in
(\ref{eqn:starfdd}).
\end{lemma}
\begin{pf}
Condition (i) ensures that each spectral measure $H^*_M$
defines a max-stable distribution with Fr{\'e}chet marginals.
Subsequently, condition (ii) ensures consistency of these
distributions (i.e., the conditions for Kolmogorov's extension theorem
are satisfied). Hence the spectral measures $H^*_M$ define a max-stable
process $X^*$ on $T$. Finally, condition (iii) ensures that
the process $X^*$ has standard Fr{\'e}chet marginals.
\end{pf}

%re5 #&#
\begin{remark}
Condition (ii) is equivalent to
%
%e13 #&#
\begin{equation}
\label{eqn:tau4} \tau^A_K = \sum
_{ J \subset M \setminus A} \tau^M_{K \cup J} \quad\quad \forall M
\in \finite(T)\setminus\{\varnothing\}, \varnothing\neq K \subset A \subset M.
\end{equation}
\end{remark}

%s3 #&#
\section{The TM process and negative definiteness of ECFs}\label
{sect:TMprocess}

For the following characterization of the set of ECFs $\Theta(T)$, we
use the fact that $\finite(T)$, the set of finite subsets of $T$, forms
a semigroup with respect to the union operation $\cup$ and with neutral
element the empty set $\varnothing$. The following notation is adopted
from \cite{molchanov_05} and \cite{bcr_84}. For a function
$f\dvtx \finite
(T) \rightarrow\RR$ and elements $K,L \in\finite(T)$, we set
\begin{eqnarray*}
(\Delta_{K}f ) (L):= f(L)-f(L\cup K).
\end{eqnarray*}

%de6 #&#
\begin{definition}[(negative definiteness/ complete alternation)] \label
{def:ND_CA}
A function $\psi\dvtx  \finite(T) \rightarrow\RR$ is called \emph{negative
definite (in the semigroup sense)} on $\finite(T)$
if for all $n \geq2$, $\{K_1,\ldots, K_n\} \subset\finite(T)$ and $\{
a_1,\ldots,a_n\} \subset\RR$ with $\sum_{j=1}^n a_j=0$
\[
\sum_{j=1}^n \sum
_{k=1}^n a_j a_k
\psi(K_j \cup K_k) \leq0.
\]
A function $\psi\dvtx  \finite(T) \rightarrow\RR$ is called \emph
{completely alternating} on $\finite(T)$
if for all $n \geq1$, $\{K_1,\ldots,  K_n\} \subset\finite(T)$ and $K
\in\finite(T)$
\[
(\Delta_{K_1}\Delta_{K_2} \cdots\Delta_{K_n} \psi )
(K) = \sum_{I \subset\{1, \ldots, n\}} (-1)^{|I|} \psi \biggl(K
\cup \bigcup_{i \in I} K_i \biggr) \leq0.
\]
\end{definition}

Because the semigroup $(\finite(T),\cup,\varnothing)$ is idempotent,
these two terms coincide. That means $\psi\dvtx  \finite(T) \rightarrow
\RR$
is \emph{completely alternating} if and only if $\psi$ is \emph
{negative definite (in the semigroup sense)}, \cf\cite{bcr_84}, 4.4.16.

%ex7 #&#
\begin{example}[(\cite{molchanov_05}, page 52)] \label{example:cap_ND}
An important example of a negative definite (completely alternating)
function on $\finite(T)$ is the \emph{capacity functional} $C\dvtx \finite
(T) \rightarrow\RR$ of a binary process $Y=\{Y_t\}_{t \in T}$ with
values in $\{0,1\}$, which is given by $C(\varnothing) = 0$ and
\[
C(A) = \PP(\exists t \in A \mbox{ such that } Y_t = 1).
\]
\end{example}

Now, we are in position to characterize the set $\Theta(T)$ of possible
ECFs on $\finite(T)$ and to define a corresponding max-linear process $X^*$.

%th8 #&#
\begin{theorem} \label{thm:ECF_ND}
\begin{enumerate}[(b)]
\item[(a)]
The function $\theta\dvtx  \finite(T) \rightarrow\RR$ is the ECF of a
simple max-stable process on $T$ if and only if the following
conditions are satisfied:
\begin{enumerate}[(iii)]
\item[(i)]$\theta$ is negative definite,
\item[(ii)]$\theta(\varnothing) = 0$,
\item[(iii)]$\theta(\{t\}) = 1$ for all $t \in T$.
\end{enumerate}
\item[(b)]  If these conditions are satisfied, the
following choice of coefficients
\begin{eqnarray*}
&&\tau^M_L := - \Delta_{\{t_1\}} \cdots
\Delta_{\{t_l\}} \theta(M \setminus L) = \sum_{I \subset L}
(-1)^{|I|+1} \theta \bigl((M \setminus L) \cup I \bigr)
\\
&&\quad \forall M \in\finite(T)\setminus\{\varnothing\}, \varnothing \neq L =
\{t_1,\ldots,t_l\} \subset M
\end{eqnarray*}
for model (\ref{eqn:starspectralmeasure}) defines a simple max-stable
process $X^*$ on $T$ which realizes $\theta$ as its own ECF $\theta^*$.
\end{enumerate}
\end{theorem}

%de9 #&#
\begin{definition}[(Tawn--Molchanov process (TM process))]
Referring to the previous work in \cite
{colestawn_96,molchanov_08,schlathertawn_02}, we will call the simple
max-stable process $X^*$ from Theorem \ref{thm:ECF_ND}\textup{(b)} \emph{Tawn--Molchanov process (TM process)} henceforth.
\end{definition}

\begin{pf*}{Proof of Theorem \ref{thm:ECF_ND}}
If $\theta$ is an ECF of a simple max-stable process $X$ on $T$, then
necessarily $\theta(\varnothing) = 0$ and $\theta(\{t\}) = 1$ for all $t
\in T$ (\cf(\ref{eqn:ECFintro})). Further, it is an application of
l'H{\^o}pitals rule that for $A \subset\finite(T)\setminus\{
\varnothing\}$
%
%e14 #&#
\begin{eqnarray}
\label{eqn:lhopital}
\nonumber
\theta(A)&=& \lim_{x \to\infty}
\frac{-\log\PP
(\bigvee_{t
\in A} X_t \leq x )}{-\log\PP(X_t \leq x)} = \lim_{x \to\infty} \frac{1-\PP (\bigvee_{t \in A} X_t \leq x
)}{1-\PP (X_t \leq x  )}
\\[-8pt]
\\[-8pt]
&=& \lim_{x \to\infty} \frac{\PP (\exists  t \in A \mbox{ such that } X_t \geq x
)}{\PP (X_t \geq x  )} = \lim_{x \to\infty}
\frac{C^{(x)}(A)}{p^{(x)}},
\nonumber
\end{eqnarray}
where $C^{(x)}$ denotes the capacity functional for the binary process
$Y_t=\mathbh{1}_{X_t \geq x}$ and $p^{(x)} = \EE Y_t = 1-\mathrm{e}^{-1/x}$.
Since negative definiteness respects scaling and pointwise limits,
negative definiteness of $\theta$ follows from Example \ref
{example:cap_ND}. This shows the necessity of (i)--(iii).

Now, let $\theta\dvtx \finite(T)\rightarrow\RR$ be a function satisfying
conditions (i)--(iii) and let the coefficients $\tau^M_L$ be given
as above. We need to check that they fulfill the (in)equalities from
Lemma \ref{lemma:tauprocess}. Indeed we have:
\begin{itemize}
\item
The inequalities $\tau^M_L = - \Delta_{\{t_1\}} \cdots\Delta_{\{t_l\}}
\theta(M \setminus L) \geq0$ follow directly from the complete
alternation of $\theta$ that is equivalent to (i).
\item From the definition of $\Delta_{\{t\}}$ we observe
\begin{eqnarray*}
\tau^{M \cup\{t\}}_{L \cup\{t\}} &=& - \Delta_{\{t\}}
\Delta_{\{t_1\}} \cdots\Delta_{\{t_l\}}\theta \bigl( \bigl(M \cup\{t\}
\bigr) \setminus \bigl(L \cup\{t\} \bigr) \bigr)
\\
&=& - \Delta_{\{t_1\}} \cdots\Delta_{\{t_l\}}\theta (M \setminus L )
+ \Delta_{\{t_1\}} \cdots\Delta_{\{t_l\}}\theta \bigl(M \cup \{ t\}
\setminus L \bigr)
\\
&=& \tau^M_L -\tau^{M \cup\{t\}}_L.
\end{eqnarray*}
\item For $t \in T$, we have $\tau^{\{t\}}_{\{t\}} = \theta(\{t\}) = 1$
because of (iii).
\end{itemize}
Thus, the coefficients $\tau^M_L$ define a simple max-stable process
$X^*$ on $T$ as given by model (\ref{eqn:starspectralmeasure}).
Finally, we compute the ECF $\theta^*$ of $X^*$ and see that it
coincides with $\theta$: For the empty set, we have $\theta
^*(\varnothing
)=0=\theta(\varnothing)$ because of (ii); otherwise we compute for $A
\subset\finite(T)\setminus\{\varnothing\}$ that
\begin{eqnarray*}
\theta^*(A) & \stackrel{{\scriptsize{(\ref{eqn:ECFintro}),(\ref
{eqn:starspectralmeasure})}}} {=}&
\sum_{\varnothing\neq L \subset A} \tau^A_L = \sum
_{\varnothing\neq L \subset A} \sum_{I \subset L}
(-1)^{|I|+1} \theta \bigl((A \setminus L) \cup I \bigr)
\\
&=& \sum_{\varnothing\neq K \subset A} \theta(K) \mathop{\sum
_{\varnothing
\neq L \subset A}}_{ A \setminus L \subset K} (-1)^{|K \cap L|+1} = \sum
_{\varnothing\neq K \subset A} \theta(K) \bigl(- (-\mathbh {1}_{K=A} )
\bigr) = \theta(A).
\end{eqnarray*}
This shows sufficiency of (i)--(iii) and part (b).\vadjust{\goodbreak}
\end{pf*}

Theorem \ref{thm:ECF_ND} is in analogy to the following standard result
for Gaussian processes (as illustrated in the sketches in the \hyperref[sect:intro]{Introduction}):
\begin{enumerate}[(b)]
\item[(a)]
A function $C\dvtx T \times T \rightarrow\RR$ is a covariance function if
and only if it is positive definite.
\item[(b)] If $C\dvtx T \times T \rightarrow\RR$ is positive definite, we may
choose a (zero mean) Gaussian process which realizes $C$ as its own
covariance function.
\end{enumerate}
Both statements are   intrinsically tied together. When proving them
by means of Kolmogorov's extension theorem, one proceeds in the same
manner as we did for Theorem \ref{thm:ECF_ND}. The necessity of
positive definiteness of covariance functions is easily derived even
for the bigger class of square-integrable processes, whilst sufficiency
can be established by showing that Gaussian processes can realize any
positive definite function as covariance function.
{In some points (such as continuity relations), this analogy will be
deepened. Other aspects (such as the exceptional role of dependency
sets in Section~\ref{sect:depset}) seem unsuitable for a direct comparison.}

%re10 #&#
\begin{remark}
In order to incorporate stationarity \wrt some group $G$ acting on $T$
(for example, $\RR^d$ acting on $\RR^d$ by translation), we just have
to add the following condition
(iv) $\theta(gA)=\theta(A)$ for all $A \in\finite(T)\setminus\{
\varnothing\}$ and for all $g \in G$.
Then the process $X^*$ will be stationary \wrt this group action.
\end{remark}

%re11 #&#
\begin{remark}
Instead of requiring the max-stable processes in Theorem \ref
{thm:ECF_ND} to have \emph{standard} Fr{\'e}chet marginals everywhere,
we can admit a different scale at different locations, that is, $\PP(X_t
\leq x)=\exp(-s_t/x)$ for a positive scaling parameter $s_t$ for $t
\in
T$. In that case Theorem \ref{thm:ECF_ND} holds true without condition
(iii) and the word ``simple''. To make sense of the ECF as in (\ref
{eqn:ECFintro}) in this case, either use a reference point $t \in T$ or
set $\log\PP(X_t\leq x)=-1/x$ in the denominator. Beware of that the
ECF $\theta$ cannot be interpreted as the number of independent
variables anymore in this case.
\end{remark}

%re12 #&#
\begin{remark}
In \cite{schlathertawn_02}, the last issue of the proof is derived for
finite sets $T$ by a Moebius inversion. The relation to the proof
therein becomes more transparent if we compute $\theta^*(A)$ for $A
\subset M$ from the coefficients $ \{\tau^M_L \}_{\varnothing
\neq L \subset M}$ for arbitrary $M \supset A$ instead of $M=A$:
%
%e15 #&#
\begin{equation}
\label{eqn:thetaA.from.tauML} \theta^*(A) \stackrel{{\scriptsize{(\ref{eqn:ECFintro}),(\ref
{eqn:starspectralmeasure})}}}
{=} \sum_{\varnothing\neq K \subset A}\tau^A_K
\stackrel{{\scriptsize{(\ref{eqn:tau4})}}} {=} \sum_{\varnothing\neq K \subset A}
\sum_{J \subset M \setminus A} \tau^M_{K \cup J} = \sum
_{L \subset M
\pmid L \cap A \neq\varnothing} \tau^M_L.
\end{equation}
\end{remark}

%s4 #&#
\section{Direct consequences of Theorem \texorpdfstring{\protect\ref
{thm:ECF_ND}}{8}}\label{sect:consequences}

Here, we collect some direct consequences of the above Theorem \ref
{thm:ECF_ND}. Therefore, note that the first part of Theorem \ref
{thm:ECF_ND} can also be expressed as (\cf(\ref{eqn:Theta_defn}))
%
%e16 #&#
\begin{equation}
\label{eqn:Theta_ND} \Theta(T)= \bigl\{\theta\dvtx \finite(T) \rightarrow[0,\infty)
\pmid %
\theta\mbox{ is negative definite, } \theta(\varnothing)=0,
\theta \bigl(\{t\} \bigr)=1 \mbox{ for } t \in T %
\bigr\}.
\end{equation}

%s4.1 #&#
\subsection{Convexity and compactness}

%co13 #&#
\begin{corollary} \label{cor:Theta_convex}
The set of ECFs $\Theta(T)$ is convex.
\end{corollary}

\begin{pf}
This can be seen directly from (\ref{eqn:Theta_ND}) since all involved
properties are compatible with convex combinations. As a constructive
argument, use the fact that the ECF of the max-combination $\alpha X
\vee(1-\alpha) Y$ of two independent simple max-stable processes $X$
and $Y$ on $T$ is the convex combination of their ECFs for $\alpha\in(0,1)$.
\end{pf}

%co14 #&#
\begin{corollary} \label{cor:Theta_compact}
The set of ECFs $\Theta(T)$ is compact \wrt the topology of pointwise
convergence.
\end{corollary}

\begin{pf} The topology of pointwise convergence on $\RR^{\finite(T)}$
is the product topology. Since $\theta(\varnothing)=0$ and $\theta
(A)\in
[1,|A|]$ for $\theta\in\Theta(T)$ and $A \in\finite(T)\setminus\{
\varnothing\}$, the set $\Theta(T)$ is a subset of the product space
\[
\{0\} \times\prod_{A \in\finite(T)\setminus\{\varnothing\}} \bigl[1,|A|\bigr],
\]
which is compact by Tychonoff's theorem. Moreover, since elements of
$\Theta(T)$ are completely characterized by finite dimensional
equalities and inequalities involving $\leq$ only (stemming from (\ref
{eqn:Theta_ND})), the set $\Theta(T)$ is closed. Hence, $\Theta(T)$
is compact.
\end{pf}

%re15 #&#
\begin{remark}
Note that even though we say ``the topology of pointwise convergence'',
the ``points'' meant here are indeed elements of $\finite(T)$, that is,
finite subsets of $T$. In particular it follows from the compactness of
$\Theta(T)$ that $\Theta(T)$ is sequentially closed. That means if
$(\theta_n)_{n \in\NN}$ is a sequence of ECFs such that $\theta_n(A)$
converges to some value $f(A)$ for each $A \in\finite(T)$, then $f$ is
an ECF.
\end{remark}

%s4.2 #&#
\subsection{Spectral representation of the TM process}

Another consequence of Theorem \ref{thm:ECF_ND} is that ECFs allow for
an \emph{integral representation} as a mixture of functions $A \mapsto
\mathbh{1}_{A \cap Q \neq\varnothing}$, where $Q$ is from the power set
of $T$. To be more precise, let us denote the power set of $T$ by
$\powerset(T)$ and consider the topology on $\powerset(T)$ that is
generated by the maps $Q \mapsto\mathbh{1}_{A \cap Q \neq\varnothing}$
for $A \in\finite(T)$ or equivalently (since $\finite(T)$ is generated
by the singletons $\{\{t\}\}_{t \in T}$) the topology on $\powerset(T)$
that is generated by the maps $Q \mapsto\mathbh{1}_{t \in Q}$ for $t
\in T$. Identifying $\powerset(T)$ with $\{0,1\}^T$, this space is also
known as \emph{Cantor cube}. As in \cite{bcr_84}, Definition 2.1.1, a
measure $\mu$ on the Borel-$\sigma$-algebra of $\powerset(T)$ \wrt this
topology will be called \emph{Radon measure} if $\mu$ is finite on
compact sets and $\mu$ is inner regular.

%co16 #&#
\begin{corollary} \label{cor:intrep}
Let $\theta\in\Theta(T)$ be an ECF. Then $\theta$ uniquely determines
a positive Radon measure $\mu$ on $\powerset(T)\setminus\{\varnothing
\}$
such that
\[
\theta(A) = \mu \bigl( \bigl\{Q \in\powerset(T) \setminus\{\varnothing\} \pmid A
\cap Q \neq\varnothing \bigr\} \bigr) = \int_{\powerset(T) \setminus\{\varnothing
\}}
\mathbh{1}_{A \cap Q \neq\varnothing} \mu(\D Q),
\]
where $\theta(\{t\})=1$ for $t \in T$.
The function $\theta$ is bounded if and only if $\mu(\powerset(T)
\setminus\{\varnothing\}) < \infty$.
\end{corollary}

\begin{pf}
Since $\theta$ is negative definite (Theorem \ref{thm:ECF_ND}) and
$\finite(T)$ is idempotent, we may apply \cite{bcr_84}, Proposition
4.4.17. It says that
$\theta$ uniquely determines a positive Radon measure $\widetilde\mu$
on $\widehat{\finite(T)}\setminus\{1\}$, where $\widehat{\finite(T)}$
denotes the dual semigroup of $\finite(T)$ (\cf\cite{bcr_84}, 4.2.1 and
4.4.16), such that $\theta(A) = \widetilde\mu(\{\rho\in
\widehat{\finite(T)} \setminus\{1\} \mid\rho(A)=0 \})$. The function
$\theta$ is bounded if and only if $\widetilde\mu(\widehat{\finite(T)}
\setminus\{1\}) < \infty$.

Now, it can be easily seen that semicharacters on $\finite(T)$ are in a
1--1 correspondence with subsets of $T$ via $\widehat{\finite(T)} \ni
\rho\rightarrow\{ t \in T \pmid\rho(\{t\})=0 \} \in\powerset(T)$
and $\powerset(T) \ni Q \rightarrow\mathbh{1}_{(\cdot) \cap Q =
\varnothing} \in\widehat{\finite(T)}$.
Here the constant function $1$ corresponds to the empty set.
Moreover, the topology considered on $\widehat{\finite(T)}$ is the
topology of pointwise convergence. Transported to $\powerset(T)$ this
is the topology generated by the maps $Q \mapsto\mathbh{1}_{A \cap Q
\neq\varnothing}$ for $A \in\finite(T)$.
Let $\mu$ denote the Radon measure $\widetilde\mu$ transported to
$\powerset(T)\setminus\{\varnothing\}$. Then the corollary follows.
\end{pf}

%re17 #&#
\begin{remark} \label{remark:tauFourier}
In case $T=M$ is finite, we have that $\powerset(M)=\finite(M)$ carries
the discrete topology and
\[
\theta(A) = \mu \bigl( \bigl\{Q \in\finite(M) \setminus\{\varnothing\} \pmid A
\cap Q \neq \varnothing \bigr\} \bigr) = \sum_{Q \in\finite(M) \setminus\{\varnothing\}}
\mu \bigl( \{Q\} \bigr) \mathbh {1}_{A \cap Q \neq\varnothing}.
\]
A comparison with (\ref{eqn:thetaA.from.tauML}) reveals that $\mu(\{
Q\}
)=\tau^M_Q$.
In this sense, the coefficients $\tau^M_Q$ of the max-linear model
(\ref
{eqn:starspectralmeasure}) can be interpreted as finite dimensional
``Fourier coefficients'' of the negative definite function $\theta$.
\end{remark}

The integral representation of the ECF $\theta$ also yields a spectral
representation for the corresponding TM process $X^*$.

%th18 #&#
\begin{theorem} \label{thm:spectralrep}
The TM process $X^*=\{X^*_t\}_{t \in T}$ with ECF $\theta$ has the
following spectral representation $(\Omega,\sigmaA,\nu,V)$ (\cf
(\ref
{eqn:spectralrepresentation})):
\begin{itemize}
\item$(\Omega,\sigmaA,\nu)$ is the measure space $(\powerset
(T),\Borel
(\powerset(T)),\mu)$ from Corollary \ref{cor:intrep},
\item$V_t(Q)=\mathbh{1}_{t \in Q}$.
\end{itemize}
\end{theorem}

\begin{pf}
We need to check that the \fdd of $X^*$ satisfy (\ref{eqn:fddspectralrep}).
The \fdd of $X^*$ are given by (\cf(\ref{eqn:starfdd}))
\begin{eqnarray*}
-\log\PP \bigl(X^*_t\leq x_t, t \in M \bigr) =\sum
_{\varnothing\neq L \subset M} \tau^M_L \bigvee
_{t \in L} \frac{1}{x_t},
\end{eqnarray*}
where the coefficients $\tau^M_L$ can be computed from the ECF $\theta$
as in Theorem \ref{thm:ECF_ND}(b) and $\theta$
satisfies the integral representation from Corollary \ref{cor:intrep},
that is,
\begin{eqnarray*}
\tau^M_L &=& \sum_{I \subset L}
(-1)^{|I|+1} \theta \bigl((M \setminus L) \cup I \bigr)
\\
 &=& \sum
_{I \subset L} (-1)^{|I|+1} \int_{\powerset(T) \setminus\{
\varnothing\}}
\mathbh{1}_{((M \setminus L) \cup I) \cap Q \neq
\varnothing
} \mu(\D Q).
\end{eqnarray*}
Using the identity\vspace*{-0.5pt}
\begin{eqnarray*}
&&\sum_{I \subset L} (-1)^{|I|+1}
\mathbh{1}_{((M \setminus L) \cup I)
\cap Q \neq\varnothing}
\\[-0.5pt]
&&\quad= \sum_{I \subset L} (-1)^{|I|+1} (
\mathbh{1}_{(M \setminus L)
\cap Q \neq\varnothing} + \mathbh{1}_{I \cap Q \neq\varnothing} - \mathbh{1}_{(M \setminus L) \cap Q \neq\varnothing}
\mathbh{1}_{I
\cap
Q \neq\varnothing} )
\\[-0.5pt]
&&\quad= 0 \cdot\mathbh{1}_{(M \setminus L) \cap Q \neq\varnothing} + (1-\mathbh{1}_{(M \setminus L) \cap Q \neq\varnothing} )
\sum_{I
\subset L} (-1)^{|I|+1}\mathbh{1}_{I \cap Q \neq\varnothing}
\\[-0.5pt]
&&\quad= \mathbh{1}_{(M \setminus L) \cap Q = \varnothing} \mathbh{1}_{L
\subset Q} =
\mathbh{1}_{L=M \cap Q},
\end{eqnarray*}
we obtain that\vspace*{-0.5pt}
\[
\tau^M_L = \int_{\powerset(T) \setminus\{\varnothing\}}
\mathbh{1}_{L=M \cap
Q} \mu(\D Q).
\]
It follows that the \fdd of $X^*$ satisfy\vspace*{-0.5pt}
\begin{eqnarray*}
-\log\PP \bigl(X^*_t\leq x_t, t \in M \bigr) &=& \int
_{\powerset(T) \setminus\{\varnothing\}} \sum_{\varnothing\neq L
\subset M}
\mathbh{1}_{L=M \cap Q} \bigvee_{t \in L}
\frac{1}{x_t} \mu (\D Q)
\\[-0.5pt]
&=& \int_{\powerset(T) \setminus\{\varnothing\}} \bigvee_{t \in M}
\frac
{\mathbh{1}_{t \in Q}}{x_t} \mu(\D Q) = \int_{\Omega} \biggl(\bigvee
_{t \in M} \frac{V_{t}(\omega
)}{x_t} \biggr) \nu(\D\omega)
\end{eqnarray*}
as desired. This finishes the proof.\vspace*{-1pt}
\end{pf}

%s4.3 #&#
\subsection{Triangle inequalities and operation of Bernstein functions}

In \cite{cooleyetalii_06}, Proposition 4, it is shown that an ECF
$\theta$ on $\finite(T)$ satisfies the following bivariate inequalities
for $r,s,t\in T$:\vspace*{-0.5pt}
\begin{eqnarray*}
\theta \bigl(\{s,t\} \bigr) & \leq&\theta \bigl(\{s,r\} \bigr)\theta \bigl(\{r,t
\} \bigr),
\\[-0.5pt]
\theta \bigl(\{s,t\} \bigr)^\alpha& \leq&\theta \bigl(\{s,r\}
\bigr)^\alpha+ \theta \bigl(\{r,t\} \bigr)^\alpha- 1,\quad\quad 0
< \alpha \leq1,
\\[-0.5pt]
\theta \bigl(\{s,t\} \bigr)^\alpha& \geq&\theta \bigl(\{s,r\}
\bigr)^\alpha+ \theta \bigl(\{r,t\} \bigr)^\alpha- 1,\quad\quad
\alpha\leq0.
\end{eqnarray*}
These inequalities have in common, that they are in fact triangle
inequalities of the form\vspace*{-0.5pt}
%
%e17 #&#
\begin{equation}
\label{eqn:BernsteinTriangle} g\circ\eta \bigl(\{s,t\} \bigr) \leq g\circ\eta \bigl(\{s,r\}
\bigr) + g\circ\eta \bigl(\{r,t\} \bigr),
\end{equation}
if we rewrite them in terms of $\eta:= \theta- 1$ and\vspace*{-0.5pt}
\begin{eqnarray*}
g(x)&=&\log(1+x),
\\[-0.5pt]
g(x)&=&(1+x)^\tau-1,\quad\quad 0 < \alpha\leq1,
\\[-0.5pt]
g(x)&=&1-(1+x)^\tau, \quad\quad\alpha\leq0.
\end{eqnarray*}
These functions $g$ have in common that they are in fact \emph
{Bernstein functions}.\vadjust{\goodbreak}

%de19 #&#
\begin{definition}[(Bernstein function)]
A function $g\dvtx [0,\infty) \rightarrow[0, \infty)$ is called a \emph
{Bernstein function} if one of the following equivalent conditions is
satisfied (\cf\cite{bcr_84}, 4.4.3 and page~141)
\begin{enumerate}[(iii)]
\item[(i)] The function $g$ is of the form
\[
g(r)= c + br + \int_0^\infty \bigl(1-
\mathrm{e}^{-\lambda r } \bigr) \nu(\D \lambda),
\]
where $c,b \geq0$ and $\nu$ is a positive Radon measure on $(0,\infty
)$ with $\int_0^\infty\frac{\lambda}{1+\lambda} \nu(\D\lambda) <
\infty$.
\item[(ii)] The function $g$ is continuous and $g \in C^\infty((0,\infty))$
with $g \geq0$ and $(-1)^n g^{(n+1)} \geq0$ for all $n\geq0$. (Here,
$g^{(n)}$ denotes the $n$th derivative of $g$.)
\item[(iii)] The function $g$ is continuous, $g \geq0$ and $g$ is negative
definite as a function on the semigroup $([0,\infty),+,0)$.
\end{enumerate}
\end{definition}

For a comprehensive treatise on Bernstein functions including a table
of examples, see \cite{schillingetalii_10}.
Bernstein functions play already an important role in the construction
of advanced Gaussian processes by generating novel covariance functions
from given ones, \cf\cite{zastavnyiporcu_11} and \cite
{porcuschilling_11}. Here, we see that they are equally useful for
generating new ECFs from given ECFs and correspondingly new
Tawn--Molchanov processes from given ones.

%co20 #&#
\begin{corollary}\label{cor:opBernstein}
Let $T$ be a set and $\theta\in\Theta(T)$ an ECF. Let $g$ be a
Bernstein function which is not constant. Then the function on $\finite(T)$
\[
A \mapsto \frac{g(\theta(A))-g(0)}{g(1)-g(0)}
\]
is again an ECF in $\Theta(T)$.
\end{corollary}

\begin{pf}
The result is immediate from Theorem \ref{thm:ECF_ND}, since Bernstein
functions operate on negative definite kernels (\cf\cite{bcr_84},
3.2.9 and 4.4.3).
\end{pf}

For instance, if $\theta$ is an ECF, then also $\log(1+\theta)/\log(2)$
or $((\theta+a)^q-a^q)/((1+a)^q-a^q)$ are ECFs for $0<q<1$ and $a \geq0$.
Finally, we show that (\ref{eqn:BernsteinTriangle}) holds true for
arbitrary Bernstein functions. In fact, the result of \cite
{cooleyetalii_06}, Proposition 4, can be generalized to the following
extent as a corollary to Theorem \ref{thm:ECF_ND}.

%co21 #&#
\begin{corollary}\label{cor:triangleineq}
Let $\theta\in\Theta(T)$ be an ECF. Set $\eta:=\theta-1$ and let $g$
be a Bernstein function. Then we have for $A,B,C \in\finite
(T)\setminus
\{ \varnothing\}$ that
\[
g\circ\eta(A \cup B) \leq g\circ\eta(C) + g\circ\eta(A \cup B) \leq g\circ\eta(A
\cup C) + g\circ\eta(C \cup B).
\]
\end{corollary}

\begin{pf}
Since $\theta$ is an ECF, it is negative definite (\cf Theorem \ref
{thm:ECF_ND}). Subtracting $1$ does not change this property. Notice
further that $\theta$ takes values in $\{0\} \cup[1,\infty)$, where
the value $0$ is only attained for the empty set $\varnothing$ (the
neutral element of $\finite(T)$). Thus, the function $\eta:=\theta-1:
\finite(T) \setminus\{\varnothing\} \rightarrow\RR$ is negative
definite and takes values only in $[0,\infty)$. Applying a Bernstein
function $g$ does not change this property (\cf\cite{bcr_84}, 3.2.9 and
4.4.3). By \cite{bcr_84}, 8.2.7, this also means that $f:=g
\circ\eta: \finite(T)\setminus\{\varnothing\} \rightarrow\RR$ is
negative definite on $\finite(T)\setminus\{\varnothing\}$. Since we have
also $f \geq0$ on $\finite(T)\setminus\{\varnothing\}$, we may derive
for $A,B,C \in\finite(T) \setminus\{\varnothing\}$
\begin{eqnarray*}
&&f(C)+f(A \cup B)-f(A \cup C)-f(C \cup B)
\\
&&\quad= \bigl(f(C)-f(A \cup C)-f(C \cup B)+f(A \cup B \cup C) \bigr)+ \bigl(f(A
\cup B)-f(A \cup B \cup C) \bigr)
\\
&&\quad=\Delta_{A}\Delta_{B}f(C) + \Delta_{C}f(A
\cup B) \leq0
\end{eqnarray*}
as desired. This finishes the proof.
\end{pf}

%s4.4 #&#
\subsection{Stochastic continuity}

In this section, we require $T$ to be a metric space. We need to define
the notion of continuity that we will use in connection with ECFs
$\theta\dvtx  \finite(T) \rightarrow[0,\infty)$. Therefore, let $f\dvtx
\finite
(T) \rightarrow\RR$ be a function on the finite subsets of $T$. Then
$f$ induces a family of functions $\{f^{(m)}\}_{m \geq0}$ where
$f^{(m)}\dvtx  T^m \rightarrow\RR$ is given by
\[
f^{(m)}(t_1,\ldots,t_m)=f \bigl(
\{t_1,\ldots,t_m\} \bigr).
\]

%de22 #&#
\begin{definition}\label{def:ECF_cts}
Let $f\dvtx  \finite(T) \rightarrow\RR$ be a function on the finite subsets
of a metric space $T$. We say that $f$ is \emph{continuous} if all
induced functions $f^{(m)}\dvtx T^m \rightarrow\RR$ are continuous for all
$m \geq0$, where $T^m$ is endowed with the product topology.
\end{definition}

%le23 #&#
\begin{lemma}\label{lemma:ECF_cty}
Let $X=\{X_t\}_{t \in T}$ be a simple max-stable process with ECF
$\theta$. Then the following implication holds:
\[
X \mbox{ is stochastically continuous} \quad\Longrightarrow \quad\theta\mbox{ is
continuous.}
\]
\end{lemma}

\begin{pf}
Stochastic continuity of $X$ means that for any $\varepsilon>
0$, for any $t \in T$ and sequence $t^{(n)} \rightarrow t$ we
have $\PP(| X_{t^{(n)}} - X_t | > \varepsilon) \rightarrow0$. From
this, we can easily derive that for any $\varepsilon>0$, any $m \in
\NN
$, any $(t_1,\ldots,t_m) \in T^m$ and a sequence $(t^{(n)}_1,\ldots
,t^{(n)}_m) \rightarrow(t_1,\ldots,t_m)$, also $\PP(\| (X_{t^{(n)}_i} -
X_{t_i})_{i=1}^m \| > \varepsilon) \rightarrow0$ for any reference
norm $\| \cdot\|$ on $\RR^m$. The latter implies the corresponding
convergence in distribution: $F_{(t^{(n)}_1,\ldots,t^{(n)}_m)}
\rightarrow F_{(t_1,\ldots,t_m)}$. Since $\log F_{(t_1, \ldots, t_m)}:
[0,\infty)^m \rightarrow\RR$ is monotone and homogeneous, we have that
for $x>0$ the point $(x,\ldots,x) \in(0,\infty)^m$ is a continuity
point of $F_{(t_1, \ldots, t_m)}$ (\cf\cite{resnick_08}, page~277).
Thus, the induced function $\theta^{(m)}$ on $T^m$ is continuous, since
$\theta^{(m)}(t_1,\ldots,t_m)= - x \log F_{(t_1, \ldots, t_m)} (x,
\ldots,
x)$. Hence, $\theta$ is continuous.
\end{pf}

Second, we prove the following upper bound that shows that stochastic
continuity of the TM process $X^*$ is indeed controlled by the
bivariate extremal coefficients.

%le24 #&#
\begin{lemma}\label{lemma:ctyestimate}
Let $X^*=\{X^*_t\}_{t \in T}$ be the TM process with ECF $\theta$. Set
$\eta:=\theta-1$.
Then we have for any $\varepsilon> 0$
\[
\PP \bigl( \bigl|X^*_s - X^*_t \bigr| > \varepsilon \bigr) \leq2
\biggl( 1 - \exp \biggl( - \frac
{\eta(\{s,t\})}{\varepsilon} \biggr) \biggr) \leq
\frac
{2}{\varepsilon} \eta \bigl(\{s,t\} \bigr).
\]
\end{lemma}

\begin{pf}
Let $\varepsilon> 0$. We will prove the statement for $2 \varepsilon$
instead of $\varepsilon$. Therefore, consider the following disjoint
events on a corresponding probability space $(\Omega, {\mathcal A},
\PP
)$ for $k=0,1,2, \ldots$\vspace*{-1pt}
\begin{eqnarray*}
A_k := \bigl\{\omega\in\Omega\pmid \bigl(X^*_s(
\omega),X^*_t(\omega) \bigr) \in(k \varepsilon, (k+2)
\varepsilon]^2 \setminus\bigl((k+1)\varepsilon, (k+2)\varepsilon\bigr]^2
\bigr\}.
\end{eqnarray*}
The disjoint union $\bigcup_{k=0}^{\infty} A_k$ is a subset of $\{
\omega\in\Omega\pmid|X^*_s(\omega) - X^*_t(\omega)| \leq
2\varepsilon\}$ and so\vspace*{-1pt}
\[
\PP \bigl( |X^*_s - X^*_t | \leq2\varepsilon \bigr)
\geq \PP \Biggl(\bigcup_{k=0}^{\infty}
A_k \Biggr) = \sum_{k=0}^\infty
\PP(A_k) = \lim_{n
\to
\infty} \sum
_{k=0}^n \PP(A_k).
\]
From (\ref{eqn:starfdd}) and Theorem \ref{thm:ECF_ND}, we see that the
bivariate distribution of the process $X^*$ is given by\vspace*{-1pt}
%
%e18 #&#
\begin{equation}
\label{eqn:starbivariate} - \log\PP \bigl(X^*_s \leq x, X^*_t
\leq y \bigr) = \frac{\eta(\{s,t\})}{x \vee
y}+\frac{1}{x \wedge y}.
\end{equation}
For further calculations, we abbreviate for $p,q \in\NN\cup\{0\}$\vspace*{-1pt}
\[
B(p,q):=\PP \bigl(X^*_s \leq p \cdot\varepsilon, X^*_t
\leq q \cdot \varepsilon \bigr).\vadjust{\goodbreak}
\]
Note that $B(p,q)=B(q,p)$ and $B(p,0)=0$. With this notation, we rearrange\vspace*{-1pt}
\[
\sum_{k=0}^n \PP(A_k) =
-B(n+1,n+1) + 2 \sum_{k=0}^n \bigl[
B(k+2,k+1)- B(k+2,k) \bigr].
\]
For the second summand, we have (\cf(\ref{eqn:starbivariate}))\vspace*{-1pt}
\begin{eqnarray*}
&& \sum_{k=0}^n \bigl[ B(k+2,k+1)-
B(k+2,k) \bigr]
\\[-1pt]
&&\quad\stackrel{{\scriptsize{(\ref{eqn:starbivariate})}}} {=} \sum
_{k=0}^n \biggl[ \exp \biggl(-\frac{1}{\varepsilon}
\biggl[\frac{\eta(\{s,t\})}{k+2} + \frac
{1}{k+1} \biggr] \biggr) - \exp \biggl(-
\frac{1}{\varepsilon} \biggl[ \frac
{\eta(\{s,t\})}{k+2} + \frac{1}{k} \biggr] \biggr)
\biggr]
\\[-1pt]
&&\quad= \sum_{k=0}^n \exp \biggl(-
\frac{1}{\varepsilon} \biggl[\frac{\eta
(\{s,t\}
)}{k+2} \biggr] \biggr) \biggl[ \exp \biggl(-
\frac
{1}{(k+1)\varepsilon
} \biggr) - \exp \biggl(-\frac{1}{k \varepsilon} \biggr) \biggr]
\\[-1pt]
&&\quad \geq\sum_{k=0}^n \exp \biggl(-
\frac{\eta(\{s,t\})}{2 \varepsilon
} \biggr) \biggl[ \exp \biggl(-\frac{1}{(k+1)\varepsilon} \biggr) - \exp
\biggl(-\frac{1}{k \varepsilon} \biggr) \biggr]
\\[-1pt]
&&\quad= \exp \biggl(-\frac{\eta(\{s,t\})}{2 \varepsilon} \biggr) \exp \biggl(-\frac{1}{(n+1) \varepsilon}
\biggr).
\end{eqnarray*}
Finally,\vspace*{-1pt}
\begin{eqnarray*}
&&\PP \bigl( \bigl|X^*_s - X^*_t \bigr| > 2\varepsilon \bigr)
\\[-1pt]
&&\quad= 1 - \PP \bigl( \bigl|X^*_s - X^*_t \bigr| \leq 2
\varepsilon \bigr) \leq1 - \lim_{n \to\infty} \sum
_{k=0}^n \PP(A_k)\vadjust{\goodbreak}
\\
&&\quad= 1 + \lim_{n \to\infty} B(n+1,n+1) - 2 \lim
_{n \to\infty} \sum_{k=0}^n
\bigl[ B(k+2,k+1)- B(k+2,k) \bigr]
\\
&&\quad\leq1 + \lim_{n \to\infty} \exp \biggl(-\frac{\eta(\{s,t\}
)+1}{(n+1)\varepsilon}
\biggr) - 2 \lim_{n \to\infty} \biggl(\exp \biggl(-\frac{\eta(\{s,t\})}{2 \varepsilon}
\biggr) \exp \biggl(-\frac{1}{(n+1)
\varepsilon} \biggr) \biggr)
\\
&&\quad= 2 - 2 \exp \biggl(-\frac{\eta(\{s,t\})}{2 \varepsilon} \biggr) \leq \frac{2}{2\varepsilon}
\eta \bigl(\{s,t\} \bigr).
\end{eqnarray*}
This finishes the proof.
\end{pf}

%th25 #&#
\begin{theorem}\label{thm:ECprocess_cty}
Let $X^*=\{X^*_t\}_{t \in T}$ be the TM process with ECF $\theta$. Then
the following statements are equivalent:
\begin{enumerate}[(iii)]
\item[(i)]$X^*$ is stochastically continuous.
\item[(ii)]$\theta$ is continuous.
\item[(iii)] The bivariate map $(s,t) \mapsto\theta(\{s,t\})$ is continuous.
\item[(iv)] The bivariate map $(s,t) \mapsto\theta(\{s,t\})$ is continuous
on the diagonal.\vadjust{\goodbreak}
\end{enumerate}
\end{theorem}

\begin{pf}
The implication $\mathrm{(i)} \Rightarrow\mathrm{(ii)}$ follows from Lemma \ref
{lemma:ECF_cty}. Clearly, continuity of $\theta$ implies continuity of
the induced function $\theta^{(2)}(s,t):=\theta(\{s,t\})$, which
implies continuity of $\theta^{(2)}$ on the diagonal. This shows the
implications $\mathrm{(ii)} \Rightarrow\mathrm{(iii)}$ and
$\mathrm{(iii)} \Rightarrow\mathrm{(iv)}$.
Finally, the implication $\mathrm{(iv)} \Rightarrow\mathrm{(i)}$ follows from Lemma \ref
{lemma:ctyestimate}, since $\eta(\{t,t\})=\theta(\{t\})-1=0$.
\end{pf}

%s5 #&#
\section{Dependency sets -- the special role of TM processes}\label
{sect:depset}

In this section, we show that the TM process $X^*$ with ECF $\theta$ is
exceptional among all max-stable processes sharing the same ECF $\theta
$ as $X^*$ in the sense that its dependency set $\DepSet^*$ (to be
introduced below) is maximal \wrt inclusion.

Therefore, recall that for a finite non-empty subset $M \subset T$ the
dependency set $\DepSet_M$ of $\{X_t\}_{t \in M}$ is the largest
compact convex set $\DepSet_M \subset[0,\infty)^M$ satisfying (\cf
(\ref{eqn:DepSetdefn}))
\[
\ell_M(x)=\sup \bigl\{ \langle x,y \rangle\pmid y \in
\DepSet_M \bigr\} \quad\quad\forall x \in[0,\infty)^M.
\]
The closed convex set $\DepSet_M$ may also be described as the
following intersection of half spaces (\cf\cite{schneider_93},
Section~1.7):
%
%e19 #&#
\begin{equation}
\label{eqn:DepSetdirect} \DepSet_M= \bigcap_{x \in S_M}
\bigl\{ y \in[0,\infty)^M \pmid\langle x,y \rangle\leq
\ell_M(x) \bigr\}.
\end{equation}

%ex26 #&#
\begin{example}[(\cite{molchanov_08}, Example 1 and Proposition 2)]
The simplest examples for dependency sets $\DepSet_M$ are the unit cube
$[0,1]^M$ corresponding to a collection of independent random variables
$\{X_t\}_{t \in M}$ and the cross-polytope $D^M:=\{x \in[0,\infty)^M
\pmid\sum_{t \in M} x_t \leq1\}$ corresponding to identical random
variables $\{X_t\}_{t \in M}$. Any dependency set $\DepSet_M$ satisfies
\[
D^M \subset\DepSet_M \subset[0,1]^M.
\]
\end{example}

\begin{figure}

\includegraphics{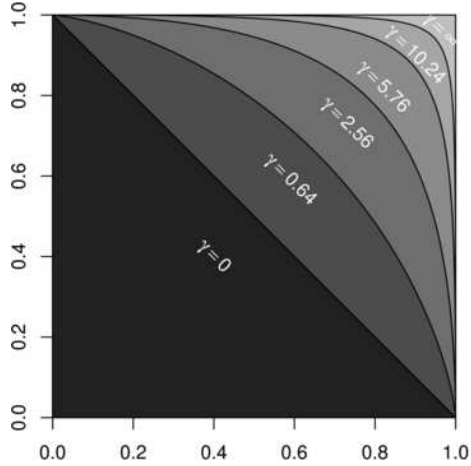}

\caption{Nested dependency sets $\DepSet^{(\gamma)}_{M}$ of
the bivariate Brown--Resnick (\resp H\"usler--Reiss)
distribution where $M=\{1,2\}$ (\cf Example \protect\ref{example:DepSetBR}).
The dependency sets grow as the parameter $\gamma$ increases.
They range between full dependence ($\gamma=0$) and independence
($\gamma=\infty$).}\label{fig:DepSetBR}
\end{figure}

%ex27 #&#
\begin{example}[(Brown--Resnick process/H\"usler--Reiss distribution)]\label{example:DepSetBR}
The \fdd of a Brown--Resnick process (\cf Example \ref{example:BR}) are
the multivariate H\"usler--Reiss distributions (\cf\cite{hueslerreiss_89}).
In the bivariate case, when $M=\{1,2\}$ consists of two points only,
the distribution function of a H\"usler--Reiss distributed random vector
$(X_1,X_2)$, standardized to unit Fr{\'e}chet marginals, is
\[
-\log\PP_{\gamma}(X_1 \leq x_1, X_2
\leq x_2) = \frac{1}{x_1} \Phi \biggl(\frac{\sqrt{\gamma}}{2} +
\frac{\log (x_2/x_1 )}{\sqrt {\gamma
}} \biggr) + \frac{1}{x_2} \Phi \biggl(\frac{\sqrt{\gamma}}{2} +
\frac{\log
(x_1/x_2 )}{\sqrt{\gamma}} \biggr)
\]
for $x_1,x_2 \geq0$. Here $\Phi$ denotes the distribution function of
the standard normal distribution and the parameter $\gamma$ is the
value of the variogram between the two points (\cf Example \ref
{example:BR}). Figure~\ref{fig:DepSetBR} illustrates, how the
corresponding dependency sets range between full dependence ($\gamma
=0$) and independence ($\gamma=\infty$).
\end{example}

%f5 #&#

In order to define a single dependency set for a simple max-stable
process comprising all multivariate dependency sets, we write
\[
\pr_M\dvtx [0,\infty)^T \rightarrow[0,
\infty)^M,\quad\quad (x_t)_{t \in T}
\mapsto(x_t)_{t \in M}
\]
for the natural projection.

%de28 #&#
\begin{definition} Let $X$ be a simple max-stable process $X=\{X_t\}_{t
\in T}$ and denote for finite $M\in\finite(T)\setminus\{\varnothing\}$
the multivariate dependency set of the random vectors $\{X_t\}_{t \in
M}$ by $\DepSet_M$. Then we define the \emph{dependency set} $\DepSet
\subset[0,\infty)^T$ of $X$ as
\[
\DepSet:= \bigcap_{M \in\finite(T)\setminus\{\varnothing\}} \pr _M^{-1}
(\DepSet_M ).
\]
\end{definition}

Analogously to (\ref{eqn:DepSetdefn}), the dependency set $\DepSet$ may
be characterized as follows.

%le29 #&#
\begin{lemma}\label{lemma:DepSetprops} The dependency set $\DepSet$ of
a simple max-stable process $X=\{X_t\}_{t \in T}$ is the largest
compact convex set $\DepSet\subset[0,\infty)^T$ satisfying
%
%e20 #&#
\begin{equation}
\label{eqn:DepSetprops} \ell_M(x)=\sup \biggl\{ \sum
_{t \in M} x_ty_t \pmid y \in\DepSet
\biggr\}\quad\quad \forall x \in[0,\infty)^M \forall \varnothing\neq M
\in \finite(T),
\end{equation}
where $\ell_M$ is the stable tail dependence function of $\{X_t\}_{t
\in M}$.
\end{lemma}

\begin{pf}
Convexity of $\DepSet$ follows from the convexity of each $\DepSet_M$
and from the linearity of the projections $\pr_M$ for $M\in\finite
(T)\setminus\{\varnothing\}$. Since $\DepSet_{\{t\}}=[0,1]$ is the unit
interval for each $t \in T$, the set $\DepSet$ is contained in the
compact space $[0,1]^T$. Moreover, $\DepSet$ is closed as the
intersection of closed sets, hence $\DepSet$ is compact.

Next, we prove that $\DepSet_M=\pr_M(\DepSet)$. By definition of
$\DepSet$ it is clear that $\pr_M(\DepSet) \subset\DepSet_M$ for $M
\in\finite(T)\setminus\{\varnothing\}$. To prove the reverse inclusion,
let $y_M$ be an element of $\DepSet_M$ and set
$V(y_M):=
\pr_M^{-1}(\{y_M\})\cap\DepSet
= \pr_M^{-1}(\{y_M\})\cap\DepSet\cap[0,1]^T$.
We need to show that $V(y_M) \neq\varnothing$.
Denoting $V(y_M,A):=\pr_M^{-1}(\{y_M\}) \cap\pr_A^{-1} (\DepSet
_A ) \cap[0,1]^T$,
we see that
\begin{eqnarray*}
V(y_M) =\bigcap_{A \in\finite(T)\setminus\{\varnothing\}}
V(y_M,A).
\end{eqnarray*}
Note that each $V(y_M,A)$ is a closed subset of the compact Hausdorff
space $[0,1]^T$. Therefore, it suffices to verify the finite
intersection property for the system of sets $\{V(y_M,A)\}_{A \in
\finite(T)\setminus\{\varnothing\}}$ in order to show that $V(y_M)\neq
\varnothing$. But this follows from the consistency of the finite
dimensional dependency sets $\{\DepSet_A\}_{A \in\finite(T)\setminus
\{
\varnothing\}}$ as follows: As \cite{molchanov_08}, Section~7, Proposition
8, essentially says, we have that if $A$ and $B$ are
non-empty finite subsets of $T$ with $A \subset B$, then $\DepSet_A$ is
the projection of $\DepSet_B$ onto the respective coordinate space. In
particular, $\pr_B^{-1}(\DepSet_B) \subset\pr_A^{-1}(\DepSet_A)$ and
$\pr_A^{-1}(\{y_A\}) \cap\pr_B^{-1}(\DepSet_B) \cap[0,1]^T \neq
\varnothing$ for $y_A \in\DepSet_A$. Now, let $A_1,\ldots,A_k$ be
non-empty finite subsets of $T$. Then
\begin{eqnarray*}
\varnothing &\neq&\pr_M^{-1} \bigl(\{y_M\}
\bigr) \cap\pr_{M \cup\bigcup_{i=1}^k A_i}^{-1} (\DepSet_{M \cup\bigcup_{i=1}^k A_i} )
\cap[0,1]^T
\\
& \subset&\pr_M^{-1} \bigl(\{y_M\} \bigr) \cap
\bigcap_{i=1}^k \pr ^{-1}_{A_i}
(\DepSet_{A_i} ) \cap[0,1]^T = \bigcap
_{i=1}^k V(y_M,A_i),
\end{eqnarray*}
as desired and we have shown that $\DepSet_M \subset\pr_M(\DepSet)$.
Both inclusions give $\DepSet_M = \pr_M(\DepSet)$.

By definition, we have $\ell_M(x)=\sup \{ \langle x,y \rangle
\pmid
y \in\DepSet_M \}$ for $x \in[0,\infty)^M$. Thus, (\ref
{eqn:DepSetprops}) follows from $\DepSet_M=\pr_M(\DepSet)$.

Finally, let $\PotentialDepSet\subset[0,\infty)^T$ be also convex
compact and satisfying (\ref{eqn:DepSetprops}) with $\DepSet$ replaced
by $\PotentialDepSet$. Then it follows immediately that $\pr
_M(\PotentialDepSet)=\DepSet_M$ for any non-empty finite subset $M
\subset T$. We conclude that $\PotentialDepSet\subset\DepSet$ by
definition of $\DepSet$. This finishes the proof.
\end{pf}

In particular, the ECF $\theta$ of a simple max-stable process $X=\{
X_t\}_{t \in T}$ can be expressed in terms of the dependency set
$\DepSet$ of $X$ as
%
%e21 #&#
\begin{equation}
\label{eqn:ECFfromDepSet} \theta(A)=\sup \biggl\{\sum_{t \in A}
x_t \pmid x \in\DepSet \biggr\}.
\end{equation}
In order to make statements about the dependency sets $\DepSet$ of
processes $X=\{X_t\}_{t \in T}$ in terms of the ECF $\theta$, we
introduce the following notation: For any non-empty finite subsets $A$
of $T$, we set the halfspace
\[
\Half_A(\theta):= \biggl\{ x \in[0,\infty)^T \pmid\sum
_{t \in A} x_t \leq\theta(A) \biggr\}
\phantom{.}
\]
that is bounded by the hyperplane
\[
\Plane_A(\theta):= \biggl\{ x \in[0,\infty)^T \pmid\sum
_{t \in A} x_t = \theta(A) \biggr\}.
\]

%le30 #&#
\begin{lemma}\label{lemma:DepSet}
Let $\DepSet$ be the dependency set of a simple max-stable process
$X=\{
X_t\}_{t \in T}$ with ECF~$\theta$. Then the following inclusion holds
\[
\DepSet\subset\bigcap_{A \in\finite(T)\setminus\{\varnothing\}} \Half _A(
\theta).
\]
On the other hand for each $A \in\finite(T)\setminus\{\varnothing\}$
there is at least one point $\mathbf{x^A}$ in the intersection
\[
\mathbf{x^A} \in\DepSet\cap\Plane_A(\theta).
\]
\end{lemma}

\begin{pf}
Let $A \in\finite(T)\setminus\{\varnothing\}$ and $x \in\DepSet$. Then
the assumption $\sum_{t \in A} x_t > \theta(A)$ contradicts $\theta
(A)=\sup\{\sum_{t \in A} x_t \pmid x \in\DepSet\} > \theta(A)$
(\cf
(\ref{eqn:ECFfromDepSet})). So $\sum_{t \in A} x_t \leq\theta(A)$.
This proves the inclusion.
Second, since $\DepSet$ is compact and the map $[0,\infty)^T \ni x
\rightarrow\sum_{t \in A} x_t$ is continuous, we know that it attains
its supremum at some $\mathbf{x^A} \in\DepSet$.
\end{pf}

%ex31 #&#
\begin{example}\label{example:ball}
We give a simple multivariate example for Lemma \ref{lemma:DepSet} (as
illustrated in Figure~\ref{fig:ballDepSet} in the introduction for the
trivariate case): The Euclidean norm $\ell_M(x)=\| x \|_2$ is a stable
tail dependence function on $[0,\infty)^M$ (\cf\cite{molchanov_08}, Example
2) and defines a simple max-stable distribution (\cf
(\ref{eqn:fddstabletaildepfn})) with ECF $\theta(A)=\sqrt{|A|}$ for $A
\subset M$, such that
\begin{eqnarray*}
\Half_A(\theta)&=& \bigl\{ x \in[0,\infty)^M \pmid
\langle x , \eins _A\rangle \leq\sqrt{|A|} \bigr\},
\\
\Plane_A(\theta)&=& \bigl\{ x \in[0,\infty)^M \pmid
\langle x , \eins _A\rangle= \sqrt{|A|} \bigr\}
\end{eqnarray*}
for $\varnothing\neq A \subset M$.
It can be easily seen that for $x \in[0,\infty)^M\setminus\{\eins
_\varnothing\}$
\[
\ell_M(x)=\| x \|_2 = \bigl\langle x, x/\| x \|_2
\bigr\rangle= \sup \bigl\{ \langle x,y \rangle\pmid y \in B^{+} \bigr\},
\]
where $B^{+}:=\{y \in[0,\infty)^M \pmid\| y \|_2 \leq1\}$ denotes
the positive part of the (Euclidean) unit ball. So, the dependency set
$\DepSet$ is clearly $B^{+}$ in this case.
Now, the planes $\Plane_A(\theta)$ are tangent to the boundary of
$B^{+}$ with common points $\mathbf{x^A}=\eins_A/\sqrt{|A|}$ for
$\varnothing\neq A \subset M$, which makes it easy to see that Lemma
\ref
{lemma:DepSet} holds true in this example. Figure~\ref{fig:ballDepSet}
shows the dependency set $\DepSet=B^+$ (left) and the intersection of
halfspaces bounded by the planes $\Plane_A(\theta)$ (right). In the
middle it is illustrated that this intersection contains $B^+$ and the
points $\mathbf{x^A}$ are marked.
\end{example}

\begin{figure}%

\includegraphics{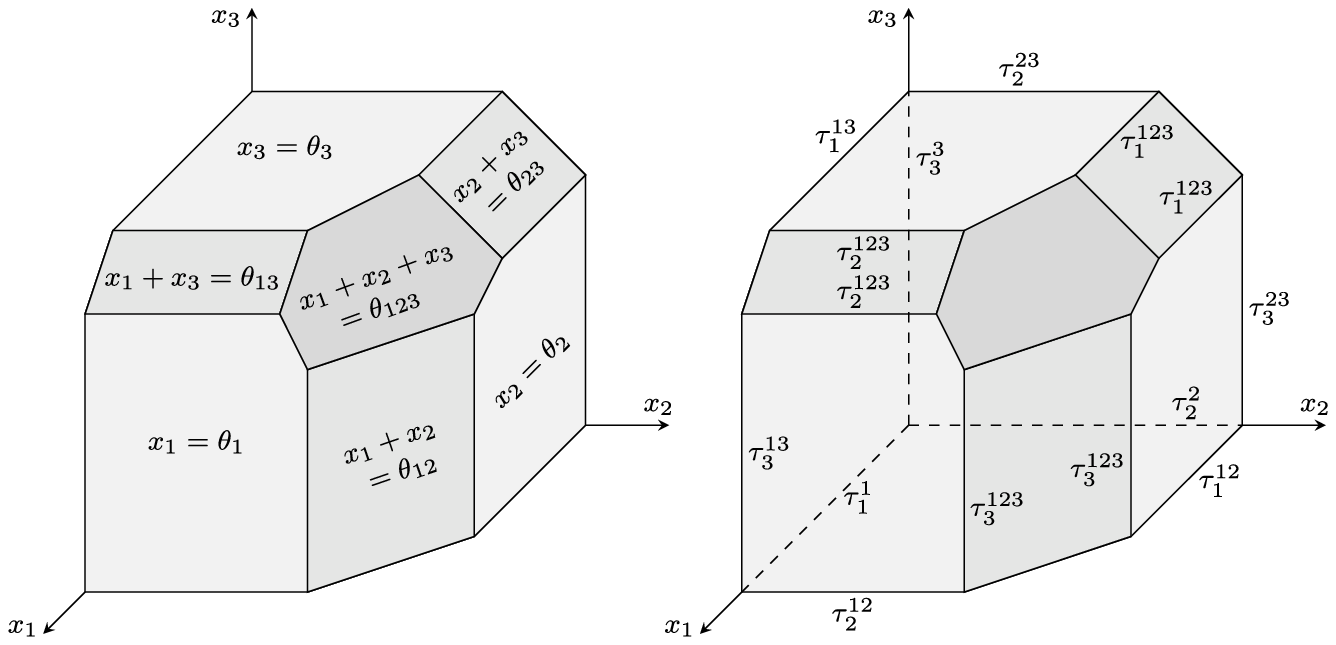}
\caption{Dependency set $\DepSet^*$ of the random
vector $\{X^*_t\}_{t \in M}$ for $M=\{1,2,3\}$. The
dependency set $\DepSet^*$ is bounded by the hyperplanes
$\Plane_A(\theta)$ that are given by the equations
$\sum_{t \in A} x_t = \theta(A)$, where $\theta$ denotes the
ECF of $X^*$. The coefficients $\tau^L_{\{t\}}$ for $L \in
\finite(M)\setminus\{\varnothing\}$ and $t \in L$ turn up as
lengths of the resulting polytope $\DepSet^*$ (\cf Theorem \protect\ref
{thm:ECF_ND}  (b) and Theorem \protect\ref
{thm:starDepSet}).} \label{fig:starDepSet}%
\end{figure}

The following theorem shows that the inclusion from Lemma \ref
{lemma:DepSet} is sharp and attained by TM processes. Figure~\ref{fig:starDepSet} illustrates the dependency set of a trivariate
distribution of a TM process.%f6 #&#

%th32 #&#
\begin{theorem}\label{thm:starDepSet}
Let $\DepSet^*$ be the dependency set of the TM process $X^*=\{X^*_t\}
_{t \in T}$ with ECF $\theta$. Then
\[
\DepSet^* = \bigcap_{A \in\finite(T)\setminus\{\varnothing\}} \Half _A(
\theta).
\]
\end{theorem}
\begin{pf}
First, we prove the theorem in the case, when $T=M$ is finite and
$\DepSet^*=\DepSet^*_M$: Therefore, write
\[
\PotentialDepSet_M := \bigcap_{\varnothing\neq A \subset M}
\Half _A(\theta) = \bigl\{ x \in[0,\infty)^M \pmid\langle
x , \eins_A\rangle \leq\theta(A) \mbox{ for all } \varnothing\neq A
\subset M \bigr\}.
\]
The inclusion $\DepSet^*_M \subset\PotentialDepSet_M$ is proven in
Lemma \ref{lemma:DepSet}. So, it remains to show the other inclusion
$\PotentialDepSet_M \subset\DepSet^*_M$. Due to (\ref
{eqn:DepSetdirect}), we have that
\[
\DepSet^*_M= \bigcap_{x \in S_M} \bigl\{ y
\in[0,\infty)^M \pmid\langle x,y \rangle\leq\ell^*_M(x)
\bigr\},
\]
where
\[
\ell^*_M(x) = \sum_{\varnothing\neq L \subset M}
\tau^M_L \bigvee_{t
\in
L}
x_{t}
\]
is the stable tail dependence function of $\{X^*_t\}_{t \in M}$, here
expressed in terms of the coefficients $\tau^M_L$ from Theorem \ref
{thm:ECF_ND} (b) (\cf(\ref{eqn:starfdd})).
Thus, it suffices to show the following implication in order to prove
$\PotentialDepSet_M \subset\DepSet^*_M$:
\[
x \in S_M \quad\mbox{and}\quad y \in\PotentialDepSet_M
\quad\Longrightarrow \quad\langle x,y \rangle\leq\ell^*_M(x).
\]
We now prove this implication: Without loss of generality, we may label
the elements of $M=\{t_1, \ldots, t_m\}$ such that $x_{t_1} \geq x_{t_2}
\geq\cdots\geq x_{t_m}$. Then we may write $x=(x_t)_{t \in M} \in S_M
\subset[0,\infty)^M$ as
\[
x = \underbrace{x_{t_m}}_{\geq0} \eins_M +
\underbrace {(x_{t_{n-1}}-x_{t_m})}_{\geq0}
\eins_{M \setminus\{t_m\}} + \cdots+ \underbrace{(x_{t_{2}}-x_{t_3})}_{\geq0}
\eins_{\{t_1,t_2\}} + \underbrace{(x_{t_{1}}-x_{t_2})}_{\geq0}
\eins_{\{t_1\}}.
\]
Taking the scalar product with $y \in\PotentialDepSet_M$, we conclude
%
%e22 #&#
\begin{eqnarray}
\label{eqn:ellcoincide} \langle x,y \rangle & \leq& x_{t_m} \theta(M) +
(x_{t_{n-1}}-x_{t_m}) \theta \bigl(M \setminus\{ t_m
\} \bigr)
\nonumber
\\
&&{} + \cdots+ (x_{t_{2}}-x_{t_3}) \theta \bigl(
\{t_1,t_2\} \bigr) + (x_{t_{1}}-x_{t_2})
\theta \bigl(\{t_1\} \bigr)
\\
& =& x_{t_m} \bigl(\theta(M)-\theta \bigl(M \setminus\{
t_m\} \bigr) \bigr) + \cdots+ x_{t_{2}} \bigl(\theta \bigl(
\{t_1,t_2\} \bigr)-\theta \bigl(\{t_1\}
\bigr) \bigr) + x_{t_{1}} \theta \bigl(\{t_1\} \bigr).
\nonumber
\end{eqnarray}
On the other hand the stable tail dependence function $\ell^*_M$ is by
this ordering of the components of $x$ given as
\[
\ell^*_M(x) = \sum_{\varnothing\neq L \subset M}
\tau^M_L \bigvee_{t
\in
L}
x_{t} = \sum_{i=1}^m
x_{t_i} \biggl( \sum_{L \subset M \pmid
t_1,\ldots
,t_{i-1} \notin L, t_{i} \in L}
\tau^M_L \biggr).
\]
From (\ref{eqn:thetaA.from.tauML}), we see that this expression
coincides with the \rhs of $(\ref{eqn:ellcoincide})$. Thus, we have our
desired inequality $\langle x,y \rangle\leq\ell^*_M(x)$. This
finishes the proof in the case, when $T=M$ is finite.

Otherwise, the definition of the dependency set $\DepSet^*$ and the
result for finite $M$ give
\[
\DepSet^* = \bigcap_{M \in\finite(T)\setminus\{\varnothing\}} \pr_M^{-1}
\bigl(\DepSet^*_M \bigr) = \bigcap_{M \in\finite(T)\setminus\{\varnothing\}}
\bigcap_{\varnothing\neq A \subset M} \pr_M^{-1}
\bigl( \Half^M_A(\theta) \bigr),
\]
where $\Half^M_A(\theta)=  \{ x \in[0,\infty)^M \pmid\sum_{t
\in
A} x_t \leq\theta(A)  \}$. Since $\pr_M^{-1}  ( \Half
^M_A(\theta)  ) = \Half_A(\theta)$ for $\varnothing\neq A
\subset
M$, the claim follows.
\end{pf}

So, if we fix the ECF $\theta$ of a simple max-stable process on $T$,
then the TM process yields a maximal dependency set $\DepSet^*$ \wrt
inclusion, that is
%
%e23 #&#
\begin{equation}
\label{eqn:DepSetinclusion} \DepSet^* = \mathop{\bigcup_{\DepSet\mbox{ {\scriptsize{dependency set}}}}}_{
\mbox{{\scriptsize{with the same ECF as}} } \DepSet^*}
\DepSet.
\end{equation}
Now, inclusion of dependency sets corresponds to stochastic ordering in
the following sense (\cf\cite{molchanov_08}, page~242): If $\DepSet
'$ and $\DepSet''$ denote the dependency sets of the simple max-stable
processes $X'$ and $X''$ respectively, then $\DepSet' \subset\DepSet
''$ implies
\[
\PP \bigl(X_{t}' \leq x_t, t \in M \bigr)
\geq\PP \bigl(X_{t}'' \leq
x_t, t \in M \bigr)\quad\quad \forall x \in[0,\infty)^M
\]
for all $M \in\finite(T)\setminus\{\varnothing\}$. This leads to the
following sharp inequality.\vadjust{\goodbreak}

%co33 #&#
\begin{corollary}\label{cor:fddinequalities}
Let $X=\{X_t\}_{t \in T}$ be a simple max-stable process with ECF
$\theta$. Let $M$ be a non-empty finite subset of $T$. Then
%
%e24 #&#
\begin{equation}
\label{eqn:DepSet_inequality} \PP(X_{t} \leq x_t, t \in M) \geq \exp
\biggl(-\sum_{\varnothing\neq L \subset M} \tau_{L}^{M}
\bigvee_{t
\in
L}\frac{1}{x_t} \biggr) \quad\quad
\forall x \in[0, \infty)^M,
\end{equation}
where the coefficients $\tau_{L}^{M}$ depend only on $\theta$ and can
be computed as in Theorem \ref{thm:ECF_ND}\textup{(b)}.
Equality holds for the TM process $X^*$.
\end{corollary}

%ex34 #&#
\begin{example}
Let us abbreviate $\eta_A:=\theta(A)-1$. In the bivariate case, the
inequality (\ref{eqn:DepSet_inequality}) reads as
\begin{eqnarray*}
\PP(X_{s} \leq x_s, X_{t} \leq
x_t) & \geq&\exp \biggl(- \biggl[\frac{\eta_{st}}{x_s \vee x_t} +
\frac{1}{x_s
\wedge x_t} \biggr] \biggr)
\\
&=& \exp \biggl(- \frac{\eta_{st}+1}{x_s \wedge x_t} \biggr) \exp \biggl(\eta _{st}
\biggl\llvert \frac{1}{x_s} - \frac{1}{x_t} \biggr\rrvert \biggr) .
\end{eqnarray*}
Indeed this inequality is much better then the trivial inequality
$\PP(X_{s} \leq x_s,  X_{t} \leq x_t) \geq\PP(X_{s} \leq x_s \wedge
x_t ,  X_{t} \leq x_s \wedge x_t)$, which can be written in the above
terms as
\[
\PP(X_{s} \leq x_s, X_{t} \leq
x_t) \geq \exp \biggl(- \frac{\eta_{st}+1}{x_s \wedge x_t} \biggr).
\]
Further note that $\eta_{st}=\theta(\{s,t\})-1$ can be interpreted as a
normalized madogram:
\[
\eta_{st}\stackrel{\scriptsize{(\ref{eqn:lhopital})}} {=}\lim
_{x \to\infty} \frac{\PP
(X_s \geq x \mbox{ or } X_t \geq x)}{\PP(X_t \geq x)}-1 = \lim_{x
\to
\infty}
\frac{\EE|\mathbh{1}_{X_s \geq x} - \mathbh{1}_{X_t \geq
x}|}{2  \EE\mathbh{1}_{X_t \geq x}}.
\]
If we additionally take into account that (\cf\cite{schlathertawn_02}, inequality
(13))
\[
\eta_{rs} \vee\eta_{st} \vee\eta_{rt} \vee(
\eta_{rs}+\eta _{st}+\eta_{rt}-1) \leq
\eta_{rst} \leq(\eta_{rs}+\eta_{st}) \wedge(
\eta_{st}+\eta_{rt}) \wedge (\eta _{rt}+
\eta_{rs}),
\]
we obtain from (\ref{eqn:DepSet_inequality}) the following (sharp)
inequality for the trivariate distribution of a simple max-stable
random vector $(X_r,X_s,X_t)$ from bivariate quantities:
\begin{eqnarray*}
&& \PP(X_{r} \leq x_r, X_{s} \leq
x_s, X_{t} \leq x_t)
\\
&& \quad\geq\exp \biggl(- \biggl[ \frac{1-\eta_{rs} \vee\eta_{st} \vee\eta_{rt}}{x_r \wedge x_s
\wedge x_t} + (a_{rst}
\wedge1) \biggl(\frac{1}{x_r \wedge x_s} + \frac{1}{x_s
\wedge x_t} + \frac{1}{x_r \wedge x_t}
\biggr)
\\
&&\quad\quad{} - \biggl( \frac{\eta_{rs}}{x_r \wedge x_s} + \frac
{\eta
_{st}}{x_s \wedge x_t} +
\frac{\eta_{rt}}{x_r \wedge x_t} \biggr) + a_{rst} \biggl(\frac{1}{x_r}+
\frac{1}{x_s}+ \frac{1}{x_t} \biggr) - \biggl(\frac{\eta_{st}}{x_r}+
\frac{\eta_{rt}}{x_s}+ \frac{\eta
_{rs}}{x_t} \biggr) \biggr] \biggr),
\end{eqnarray*}
where $a_{rst}:=(\eta_{rs}+\eta_{st}) \wedge(\eta_{rs}+\eta_{rt})
\wedge(\eta_{st}+\eta_{rt})$.
\end{example}

Thus, if one can handle the ECF of a max-stable process, sharp lower
bounds for its \fdd are available.\vadjust{\goodbreak} However, beware that higher variate
cases of these inequalities will be numerically unstable.

%re35 #&#
\begin{remark}
It is an open problem and it would be interesting to know whether there
exist also minimal dependency sets in the sense of (\ref
{eqn:DepSetinclusion}) and if they would help to better understand the
classification of all dependency structures. In view of Lemma \ref
{lemma:DepSet} and Theorem \ref{thm:starDepSet} a very naive idea would
be to take one point from each of the sets $\DepSet^* \cap\Plane_A$
where $A \in\finite(T)\setminus\{\varnothing\}$ and then to take the
convex hull with 0 included. However, this fails to be a dependency set
in dimensions $|T| \geq3$, since it is not even a zonoid, which would
be necessary (\cf\cite{molchanov_08}).
\end{remark}

\section*{Acknowledgements}
We would like to thank Ilya Molchanov for an inspiring
discussion and Zakhar Kabluchko for pointing us to the
Cantor cube. We are grateful to two unknown referees for
their valuable hints and comments that helped to significantly
improve the paper. Financial support for K. Strokorb by the German
Research Foundation DFG through the Research Training Group 1023 and
for M. Schlather by Volkswagen
Stiftung within the ``WEX-MOP'' project is gratefully acknowledged.

%suskaldyti doi

% imsref loaded by arune.pranskunaite, 2014-03-05 08:44:53
%

\printhistory

\end{document}